\documentclass[%
 reprint,
 amsmath,amssymb,
 aps,pre
]{revtex4-2}

\usepackage{graphicx}
\usepackage{dcolumn}
\usepackage{bm}
\usepackage[ruled]{algorithm2e}
\usepackage{tikz}

\DeclareUnicodeCharacter{2009}{\,} 
\begin{document}

\preprint{APS/123-QED}

\title{Tensor approximation of functional differential equations}

\author{Abram Rodgers}
\email{abram.k.rodgers@nasa.gov}
\author{Daniele Venturi}%
\email{venturi@ucsc.edu}
 
\affiliation{%
$^*$Advanced Supercomputing Division, NASA Ames Research Center,\\
$^\dagger$Department of Applied Mathematics, University of California Santa Cruz%
}%

\date{\today}

\begin{abstract}
Functional Differential Equations (FDEs) play a fundamental role in many areas of mathematical physics, including fluid dynamics (Hopf characteristic functional equation), quantum field theory (Schwinger-Dyson equation), and statistical physics. Despite their significance, computing solutions to FDEs remains a longstanding challenge in mathematical physics. In this paper we address this challenge by introducing new approximation theory and high-performance computational algorithms designed for solving FDEs on tensor manifolds. Our approach involves approximating FDEs using high-dimensional partial differential equations (PDEs), and then solving such high-dimensional PDEs on a low-rank tensor manifold leveraging high-performance parallel tensor algorithms. The effectiveness of the proposed approach is demonstrated through its application to the Burgers-Hopf FDE, which governs the characteristic functional of the stochastic solution to the Burgers equation evolving from a random initial state. 
\end{abstract}


\maketitle

\section{Introduction}

Functional differential equations (FDEs) are 
equations involving operators (nonlinear functionals), 
and derivatives/integrals of such operators with respect 
functions (functional derivatives) and other independent 
variables such as space and time \cite{venturi2018numerical,VenturiSpectral}. 
FDEs arise naturally in many different branches of physics.
A classical example in fluid dynamics is
the Hopf-Navier-Stokes equation \cite{Hopf,Monin2,Ohkitani,Rosen_1971} 
\begin{align}
\frac{\partial \Phi([ \bm \theta],t)}{\partial t}=
\sum_{k=1}^3
\int_{\Omega}\theta_k( \bm x)
&\left  (
i \sum_{j=1}^3\frac{\partial }{\partial x_j}
\frac{\delta^2 \Phi([ \bm \theta],t)}{\delta \theta_k( \bm x)\delta\theta_j(\bm x)}
+ \right.
\nonumber\\
&
\nu \nabla^2\frac{\delta \Phi([ \bm \theta],t)}{\delta \theta_k(\bm  x)}
\Biggl) d \bm x, 
\label{hopfns}
\end{align}
which governs the temporal evolution of the characteristic functional 
\begin{equation}
\Phi([ \bm \theta],t)=\mathbb{E}\left\{\exp \left(i\int_{\Omega}
 \bm u (\bm x,t)\cdot \bm \theta(\bm x)d \bm x\right)\right\}.
\label{HcF}
\end{equation}
Here, $\Omega\subseteq\mathbb{R}^3$ is the spatial 
domain, $\bm u(\bm x,t)$ represents a stochastic 
solution of the Navier-Stokes equations 
\begin{equation}
\begin{cases}
\displaystyle\frac{\partial \bm u}{\partial t}+(\bm u\cdot \nabla) \bm u = - \nabla p +\nu\nabla^2 \bm u,  \vspace{0.2cm}\\ 
\nabla \cdot \bm u= 0,
\label{NS}
\end{cases}
\end{equation}
corresponding to a random initial velocity distribution 
\cite{VakhaniaBOOK,Foias1}, $\mathbb{E}\{\cdot\}$ 
is the expectation over the probability measure of 
such random velocity distribution, and 
$\bm \theta(\bm x)$ is a divergence-free test 
function. Equation \eqref{hopfns} involves 
derivatives of the functional $\Phi([\bm \theta],t)$ with respect to the 
functions $\theta_i(\bm x)$, i.e., functional 
derivatives $\delta/\delta \theta_i(\bm x)$ \cite{VenturiSpectral}, 
and derivatives with respect to the independent variables $x_j$ and $t$. 

As is well-known, the characteristic functional \eqref{HcF} 
encodes the {full statistical information} of the random field 
$\bm u(\bm x,t)$ that solves the Navier-Stokes equations \eqref{NS}, 
including multi-point statistical moments, cumulants, 
and multi-point joint probability density functions.
For instance, by taking functional derivatives of $\Phi$ 
we can immediately express the mean and auto-correlation
function of $\bm u(\bm x,t)$ as
\begin{align}
\mathbb{E}\left\{u_k(\bm x,t)\right\}&=\frac{1}{i}
\left.\frac{\delta\Phi([\bm \theta],t)}
{\delta \theta_k( \bm x)}\right|_{\bm \theta=\bm 0},
\nonumber\\
\mathbb{E}\left\{u_k(\bm x,t)u_j(\bm y,t)\right\} &= 
\frac{1}{i^2}\left.
\frac{\delta^2 \Phi([ \bm \theta],t)}{\delta 
\theta_k( \bm x)\delta\theta_j(\bm y)}\right|_{\bm \theta=\bm 0}.
\end{align}
Similarly, by evaluating  \eqref{HcF} at 
$\theta_i(\bm x) = a_i \delta(\bm x-\bm y)$ yields 
the one-point characteristic function of the solution $\bm u(\bm y,t)$ 
at an arbitrary spatial location $\bm y\in \Omega$ as 
\begin{align}
\phi(\bm a,\bm y, t)&= \mathbb{E}\left\{\exp \left(i\int_{\Omega}
 \bm u (\bm x,t)\cdot \bm a \delta (\bm x-\bm y)d \bm x\right)\right\}\nonumber \\
&=\mathbb{E}\left\{e^{i\bm a\cdot \bm u(\bm y,t)}\right\}.
\end{align}
The inverse Fourier transform (in the sense 
of distributions) of $\phi(\bm a,\bm y, t)$
with respect to the variable $\bm a$ is the 
probability density function PDF of the solution 
$\bm u(\bm y,t)$ at the spatial location $\bm y\in \Omega$ 
and time $t$. Any other statistical property of $\bm u(\bm x,t)$ can be 
derived from the characteristic functional \eqref{HcF}.
For this reason, the FDE \eqref{hopfns} was deemed 
by Monin and Yaglom (\cite[Ch. 10]{Monin2}) 
to be {\em ``the most compact formulation of the general 
turbulence problem''}, which is the problem of 
determining the statistical properties of the velocity 
and the pressure fields of the Navier-Stokes equations 
given statistical information on the initial state. 

Another well-known classical example of functional differential 
equation is the Schwinger-Dyson (SD) equation of quantum 
field theory \cite{Peskin,Justin}. The SD equation describes 
the dynamics of the generating functional 
of the Green functions of a quantum field theory, 
allowing us to propagate quantum field interactions 
in a perturbation setting (e.g., with Feynman 
diagrams), or in a strong coupling regime. 
The usage of FDEs saw rapid growth in the 1970s, 
driven by the realization that techniques originally 
developed for quantum field theory by Dyson, Feynman, 
and Schwinger could be, at least formally, extended 
to other branches of mathematical physics. A pivotal contribution 
to this evolution was the groundbreaking work by Martin, Siggia, 
and Rose \cite{Martin}. This work served as a landmark in the field, 
showcasing the possibility of applying quantum field theoretic 
methods, such as functional integrals and diagrammatic 
expansions \cite{Phythian, Jensen, Phythian1, Jouvet}, to 
problems in classical statistical physics. 

More recently, FDEs have appeared in both mean-field games 
and mean-field optimal control \cite{Ruthotto2020,Weinan2019}. 
Mean-field games are optimization problems involving 
a large (potentially infinite) number of interacting 
players. In some cases, it is possible to reformulate such 
optimization problems in terms of a nonlinear 
Hamilton-Jacobi-Bellman FDE in a Wasserstein space 
of probability measures.
The standard form of such equation is (see \cite[p. 1]{Osher2019_1})
\begin{align}
\frac{\partial F([\rho],t)}{\partial t}&+
\mathcal{W}\left([\rho],\left[\frac{\delta F([\rho],t)}{\delta \rho(\bm x)}\right]\right)=0,\nonumber \\
F([\rho],0)&=F_0([\rho]),
\label{FDE11}
\end{align}
where $\rho(\bm x)$ is a $d$-dimensional probability density 
function supported on $R\subseteq \mathbb{R}^d$, 
$\delta F /\delta \rho(\bm x)$ is the first-order functional 
derivative of $F$ relative to $\rho(\bm x)$, and $\mathcal{W}$ 
is the Hamilton functional
\begin{align}
\mathcal{W}\left([\rho],\left[\frac{\delta F([\rho],t)}{\delta \rho(\bm x)}\right]\right) = 
\int_{R} &{\Psi}\left(\bm x,\nabla \frac{\delta F([\rho],t)}
{\delta \rho(\bm x)} \right)\rho(\bm x)d\bm x
\nonumber \\
&+\mathcal{J}([\rho]).
\end{align}
Here, $\Psi$ is the Hamilton's function and 
$\mathcal{J}([\rho])$ is an interaction potential.
Mean-field theory is also useful in optimal 
feedback control of nonlinear stochastic 
dynamical systems and in deep learning. 
For instance, E and collaborators \cite{Weinan2019} 
established the mathematical groundwork for the population risk 
minimization problem in deep learning by framing it as a 
mean-field optimal control problem. This yields 
yields a generalized version of the Hamilton-Jacobi-Bellman 
equation in a Wasserstein space of probability measures, 
which is a nonlinear FDE of the form 
\eqref{FDE11} -- e.g.,  equation (20) in \cite{Weinan2019} 
and equation (1.1) in \cite{Ganbo2021}. \\

In this paper we develop approximation theory and high-performance 
computational algorithms designed for solving FDEs on tensor manifolds.
Our approach involves first approximating FDEs in terms of 
high-dimensional partial differential equations (PDEs), and 
then computing the solution to such high-dimensional 
PDEs on a low-rank tensor manifold leveraging high-performance 
(parallel) tensor algorithms. To this end, we will build upon 
our recent work on tensors,
in particular step-truncation tensor methods
\cite{rodgers2020step-truncation,rodgers2023implicit}, 
and demonstrate convergence to functional approximations of FDEs.

This paper is organized as follows. In Section \ref{sec:exist_and_unique} 
we briefly review necessary and sufficient conditions for 
existence and uniqueness of solutions to FDEs. In Section 
\ref{sec:FTE_to_hyperPDEs} we discuss functional 
approximation of FDEs in terms of high-dimensional PDEs. 
In Section \ref{sec:approximationPDEtensor} we briefly review 
numerical tensor methods approximate the solution of 
high-dimensional PDEs on tensor manifolds.
In Section \ref{sec:approximationPDEtensor} we demonstrate 
the new functional approximation methods for FDEs 
to the Burgers-Hopf equation. The main findings
are summarized in Section \ref{sec:conclusion}.
We also include Appendix \ref{app:A} 
in which we describe the high-performance 
(parallel) tensor algorithms we developed to solve 
high-dimensional PDEs corresponding to FDEs 
on tensor manifolds.

\section{Existence and uniqueness of solutions to FDEs}
\label{sec:exist_and_unique}

Let us briefly describe what the problem is here, 
and our intended course of action. To this end, consider 
the Hopf-Navier-Stokes functional differential 
equation \eqref{hopfns}. Back in 1972 Monin and Yaglom 
stated in \cite[p. 773]{Monin2} that:

\vspace{0.2cm}
\noindent
{\em ``When we tried to develop a complete 
statistical description of turbulence with the 
aid of the Hopf equation for the characteristic functional 
we found that no general mathematical formalism for solving 
linear equations in functional derivatives was available. 
There are also no rigorous theorems on the existence 
and uniqueness of the solution to such equations.''}

\vspace{0.2cm}
\noindent
Since then, there were advancements 
in the theory of existence and uniqueness of 
solutions to the FDE \eqref{hopfns}. In particular, 
it was found the {existence} of a solution to \eqref{hopfns} 
is strongly related to the existence 
of a Hopf-Leray weak solution for the Cauchy 
problem of the Navier--Stokes equations, which is a 
well-established result \cite{Berselli,Prodi,Leray,Fabes} 
(see also \cite[Theorem 3.1]{Temam}).
Specifically, it can be proved that if there exists a Hopf-Leray solution 
to the Navier-Stokes equations then there exists a statistical 
solution to the Hopf-Navier-Stokes equation \eqref{hopfns} 
(see \cite[\S 5 in Ch. 4]{Vishik_book_1988}). 
A similar results was obtained by Foias 
in \cite[p. 254]{Foias1} for initial velocity 
conditions described by probability measures
with finite second moments.

Regarding uniqueness of the solution to the Hopf-Navier-Stokes equation
\eqref{hopfns}, it turns out that this is again connected to the uniqueness 
of Hopf-Leray weak solutions to the Navier-Stokes equation. 
Indeed,  Vishik and Fursikov proved in \cite[Theorem 6.1]{Vishik_book_1988} 
that if there exists a unique Hopf-Leray solution to the 
Navier-Stokes equation then any space-time statistical solution 
of the Hopf-Navier-Stokes equation is uniquely determined by 
the initial probability measure (see also \cite[pp. 323-344]{Foias1}).
In two spatial dimensions it is well-known that there exists 
a {unique} global Hopf-Leray solution, and therefore 
a unique statistical solution. The proof was formalized 
by Gishlarkaev in \cite{Gishlarkaev} for the 
Hopf-Navier-Stokes equation using techniques 
developed in \cite{Vishik_book_1988}.
The uniqueness of the Hopf-Leray solution for the
three-dimensional Navier-Stokes equations is still 
an open problem, and a rather active area of research. For example, it can 
be shown that the Hopf-Leray solution 
is unique within at least a short-time integration period 
if we choose the initial condition in non-conventional 
function spaces \cite{Barker2018,Barker2020,Leray,Prodi,Fabes}. 
Also, the Hopf-Leray solution is unique and global in time if we choose  
the initial data small enough \cite{Kato1984,Fabes}. 
Recently, Galdi proved in \cite{Galdi} that very 
weak solutions to the Cauchy problem for the 
Navier-Stokes equations must be Hopf-Leray 
solutions if their initial data are solenoidal 
with finite kinetic energy. 
A smooth (mild or strong) solution of the 
Navier-Stokes equations is of course also a Hopf-Leray (weak)
solution. More precisely, it can be shown that a mild solution 
to Navier-Stokes in $H^s$ agrees with a weak Hopf-Leray 
solution almost everywhere in $[0,T^*)$. This result, which 
belongs to  the so-called weak-strong uniqueness 
methods \cite{Dubois}, implies that weak Hopf-Leray solutions 
are unique as long as a mild (or a strong) solution exists. 
In particular, Leray-Hopf solutions are unique whenever 
they are regular enough to be strong solutions. 
Once the strong solution is integrated past 
$t\sim T^*$ (maximum Cauchy development), the is no 
further guarantee of uniqueness \cite{FujitaKato,Fabes}. 
We also recall that if the energy inequality that characterizes 
Hopf-Leray weak solutions is dropped then it is possible to construct 
multiple weak solutions of the Navier-Stokes equations 
with {finite} kinetic energy \cite{Buckmaster2019}. 
On the other hand, weak solutions obeying the energy 
inequality are necessarily unique and smooth 
if they lie in $L^p([0,T])\otimes L^q(\Omega)$ with $2/p+3/q=1$ 
and $q>3$. These conditions are 
known as Prodi-Serrin-Ladyzhenskaya conditions, 
and they have been recently generalized in various ways
(see \cite{Prodi,Kang,Serrin}).

\section{Approximation of FDEs by high-dimensional PDEs}
\label{sec:FTE_to_hyperPDEs}

Consider a well-posed initial/boundary value problem 
for an FDE, i.e., a problem that admits a unique solution in 
some function space as discussed in Section \ref{sec:exist_and_unique}. 
How do we compute the solution? This is a longstanding problem 
in mathematical physics and a rather new research area of 
computational mathematics. 

The FDEs \eqref{hopfns} and \eqref{FDE11} can be seen 
as PDEs in an infinite number of independent variables. 
Such infinite-dimensional PDEs may be approximated 
by PDEs in a finite (though possibly extremely large) number of 
variables using the functional analytic methods described in \cite{VenturiSpectral,venturi2018numerical,Ganbo2021,Osher2019_1}. 
For instance, we have shown in \cite{venturi2018numerical,VenturiSpectral}
that equation \eqref{hopfns} can be approximated by the high-dimensional linear PDE 
\begin{equation}
\frac{\partial \phi}{\partial t}= \sum_{p=1}^n a_p\left(
\nu \sum_{k=1}^n B_{pk}\frac{\partial \phi}{\partial a_k}+i \sum_{j,k=1}^n A_{pjk} 
\frac{\partial^2\phi}{\partial a_k\partial a_j}\right),
\label{hyper1}
\end{equation}
where 
\begin{align}
a_p=\int_{\Omega} \bm \theta\cdot \bm \Gamma_p d\bm x, &\qquad 
B_{pk}=\int_\Omega {\bm \Gamma}_p\cdot 
\nabla^2{\bm \Gamma}_k d \bm x, \nonumber\\ 
A_{pjk}&=\int_\Omega {\bm \Gamma}_p\cdot 
\left[\left({\bm \Gamma}_k\cdot \nabla\right){\bm \Gamma}_j \right]d \bm x,
\label{coefficients}
\end{align}
$\{ \bm \Gamma_1(\bm x),\ldots, \bm \Gamma_n(\bm x)\}$ 
is a divergence-free basis, and $\phi(a_1,\ldots,a_n,t)$ is 
the characteristic function approximating the characteristic functional 
$\Phi([\bm \theta],t)$ in the following sense:  
\begin{equation}
\phi(a_1,\ldots,a_n,t) =\Phi([P_n \bm \theta],t), 
\label{sense}
\end{equation}
where $P_n$ is an orthogonal projection 
onto $\{ \bm \Gamma_1,\ldots, \bm \Gamma_n\}$.
The divergence-free basis 
$\{ \bm \Gamma_1,\ldots, \bm \Gamma_n\}$ can be 
constructed in terms of div-free wavelets \cite{Deriaz1,Deriaz}, 
radial basis functions \cite{Fuselier}, trigonometric 
polynomials \cite{Landriani}, or eigenvalue problems 
for the Stokes operator with appropriate 
boundary conditions \cite{Berselli}.
Approximating the FDE \eqref{hopfns} using the PDE \eqref{hyper1} 
bears a resemblance, in a certain aspect, to what is known as the 
``method of lines'' for discretizing a nonlinear partial differential 
equation into a high-dimensional system of nonlinear ordinary 
differential equations (ODEs). However, in the context of FDEs, 
the ``lines'' are determined by a hyper-dimensional PDE instead 
of a system of ODEs.

In \cite[Theorem 7.1]{VenturiSpectral} we proved that if 
we choose the test function $\bm \theta(\bm x)$ in a 
Sobolev sphere $\mathcal{S}(\Omega)$ of $H^1(\Omega)$-div 
functions  (which is a compact subset of $L^2(\Omega)$-div) 
then $\phi$ converges uniformly to the characteristic 
functional $\Phi$ as the number of independent variables $n$ 
goes to infinity, i.e.,
\begin{equation}
\lim_{n\rightarrow \infty}\sup_{\substack{\bm \theta\in \mathcal{S}(\Omega)\\
t\in [0,T]}}\left|\phi(a_1,\ldots, a_n,t)-\Phi([\bm \theta],t)\right|=0.
\label{conv1}
\end{equation}
Moreover, convergence of the high-dimensional PDE to the FDE 
can be exponential. The result \eqref{conv1} leverages 
the Trotter-Kato approximation theorem \cite[p. 210]{engel1999one} 
for abstract evolution equations in Banach spaces \cite{Guidetti2004}.

Gangbo {\em et al.} \cite{Ganbo2021} proved a similar result for 
the nonlinear functional HJB equation \eqref{FDE11}. Specifically, 
by evaluating the functional equation \eqref{FDE11} 
on the Wasserstein space $\mathcal{P}_2\left(\mathbb{R}^d\right)$ 
(space of probability densities in $\mathbb{R}^d$ 
with finite second-order moment) with elements 
\begin{equation}
\rho_n(\bm x)=\frac{1}{n}\sum_{k=1}^n \delta (\bm x-\bm x_k)\qquad \bm x,\bm x_i \in \mathbb{R}^d
\label{empmeasure}
\end{equation}
then we can write the FDE \eqref{FDE11} as a 
hyper-dimensional nonlinear PDE
\begin{equation}
\frac{\partial f}{\partial t}+
\frac{1}{n}\sum_{k=1}^n \Psi\left(\bm x_i,\frac{\partial f}
{\partial \bm x_k}\right)+j(\bm x_1,\ldots, \bm x_n)=0,
\label{hyper2}
\end{equation}
where $f(\bm x_1,\ldots,\bm x_n,t)=F([\rho_n],t)$ 
(see \cite[Eq.(1.2)]{Ganbo2021}).
Moreover, if we choose the empirical measure \eqref{empmeasure} 
in a bounded subset ${B}$ of the Wasserstein space 
$\mathcal{P}_2\left(\mathbb{R}^d\right)$ 
then $f(\bm x_1,\ldots,\bm x_n,t)$ 
converges uniformly to the (unique)
solution of \eqref{FDE11} as $n$ goes to 
infinity (\cite[Theorem 1.2]{Ganbo2021}), i.e., 
\begin{equation}
\lim_{n\rightarrow \infty}\sup_{\substack{\rho_n\in {B}
\\ t\in [0,T]}}\left|f(\bm x_1,\ldots, \bm x_n,t)-F([\rho_n],t)
\right|=0.
\label{conv2}
\end{equation}
The boundedness assumption of $B$ is essential for convergence. 
Similarly, the compactness assumption 
of the function space $S$ in \eqref{conv1} is essential for convergence 
(see \cite[Theorems 7.1 and 8.3]{VenturiSpectral}).

To solve high-dimensional PDEs such as \eqref{hyper1} and 
\eqref{hyper2} numerically, we need a computational paradigm 
for linear PDEs that can handle tens, hundreds, or thousands 
of independent variables efficiently. Such computational 
paradigm is discussed in the next section.

\section{Tensor approximation of high-dimensional PDEs}
\label{sec:approximationPDEtensor}

\begin{figure*}[t]
\centerline{
\includegraphics[scale=0.50]{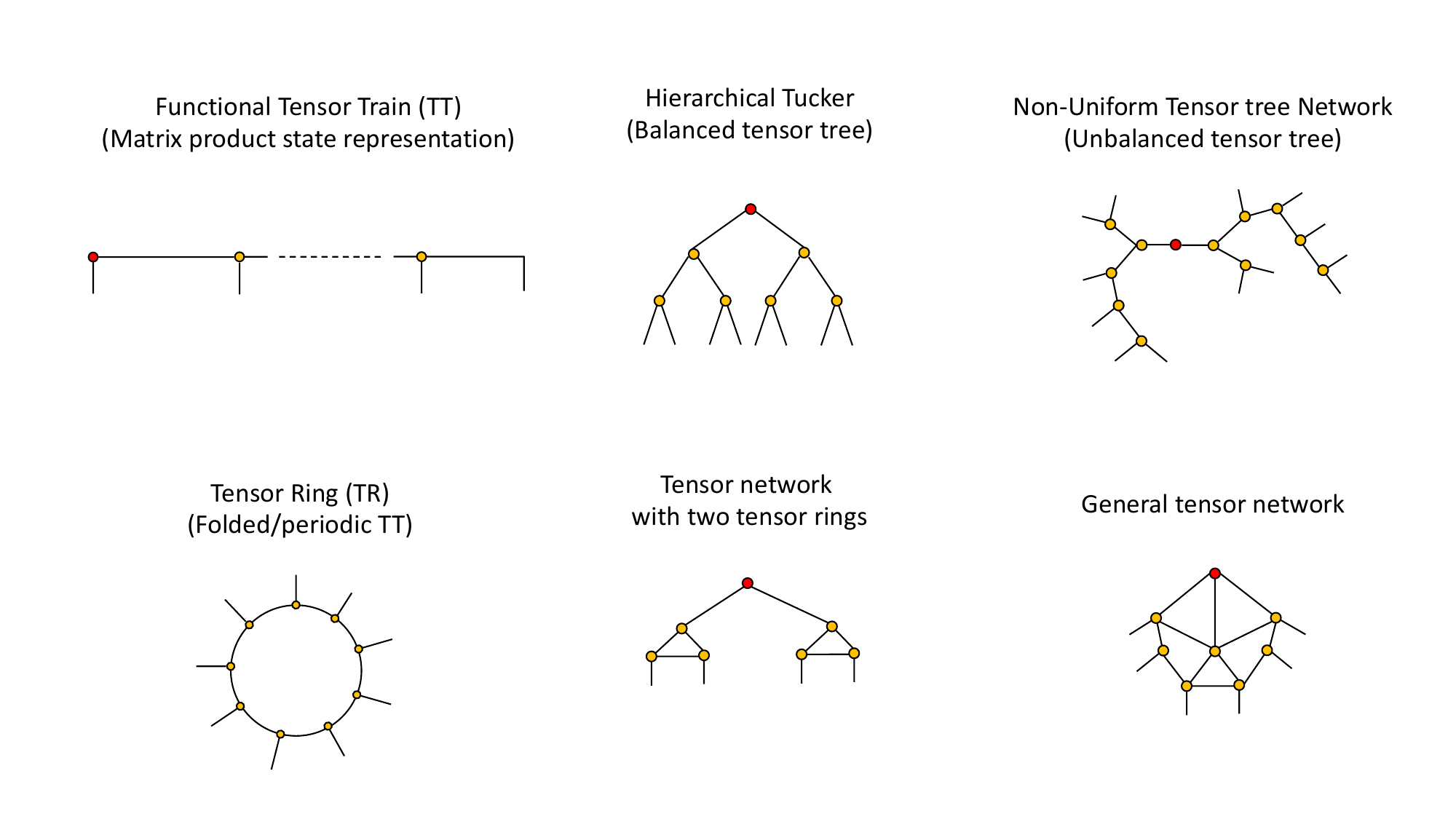}
}
\vspace{-0.5cm}
\caption{Examples of tensor networks. The vertices represent tensor 
modes (functions) used in the decomposition.  The edges connecting to vertices 
represent summation over an index (contraction) between two modes. 
The free edges represent input variables $x_j$ ($j=1,\ldots,n$).}
\label{fig:TN}
\end{figure*}

We have seen in Section \ref{sec:FTE_to_hyperPDEs} that 
FDEs can be approximated by high-dimensional PDEs of the form
\begin{align}
\frac{\partial f({\bm x},t) }{\partial t} = 
{\cal G}\left(f({\bm x},t),{\bm x}\right), \qquad 
f({\bm x},0) = f_0({\bm x}),
\label{nonlinear-ibvp} 
\end{align}
where $f:  \Omega \times [0,T] \to\mathbb{R}$ 
is a $n$-dimensional time-dependent scalar field 
defined on the domain $\Omega\subseteq \mathbb{R}^n$, 
$T$ is the period of integration, 
and $\cal G$ is a nonlinear operator which may depend 
on the variables ${\bm x}=(x_1,\ldots,x_n)\in \Omega$, and may 
incorporate boundary conditions.
For simplicity, we assume that the domain $\Omega$ 
is a Cartesian product of $n$ one-dimensional 
domains $\Omega_i$
\begin{equation}
\Omega=\Omega_1\times\cdots\times \Omega_n,
\end{equation}
and that $f$ is an element of a Hilbert 
space $H(\Omega;[0,T])$. In these hypotheses, we can leverage 
the isomorphism  $H(\Omega;[0,T])\simeq H([0,T])\otimes H(\Omega_1)\otimes \cdots \otimes H(\Omega_n)$ and represent the 
solution of \eqref{nonlinear-ibvp} as 
\begin{equation}
f(\bm x,t) \approx \sum_{i_1=1}^{m_1}\cdots
\sum_{i_n=1}^{m_n}f_{i_1\ldots i_n}(t) 
h_{i_1}(x_1)\cdots h_{i_1}(x_1),
\label{tt}
\end{equation}
where $h_{i_j}(x_j)$ are one-dimensional 
orthonormal basis functions of $H(\Omega_i)$. 
Substituting \eqref{tt} into \eqref{nonlinear-ibvp} and  
projecting onto an appropriate finite-dimensional subspace 
of $H(\Omega)$ yields the semi-discrete form 
\begin{equation}
\label{eqn:ode}
\frac{d{\bm f}}{d t}
=
{\bm G}({\bm f}),\qquad {\bm f}(0) = {\bm f}_0
\end{equation}
where 
${\bm f}:[0,T]\rightarrow
{\mathbb R}^{m_1\times m_2\times \dots \times m_n}$ is 
a multivariate array with coefficients 
$f_{i_1\ldots i_d}(t)$, and $\bm G$ 
is the finite-dimensional representation of the 
nonlinear operator $\cal G$.
The number of degrees of freedom associated 
with the solution to the Cauchy problem \eqref{eqn:ode} 
is $N_{\text{dof}}=m_1 m_2 \cdots m_n$ at each 
time $t\geq 0$, which can 
be extremely large even for 
moderately small dimension $n$. For instance, 
the solution of the Boltzmann-BGK equation
on a six-dimensional ($n=6$) flat torus \cite{BoltzmannBGK2020,dimarco2014}  
with $m_i=128$ basis functions in each position and momentum variable  
yields $N_{\text{dof}}=128^6=4398046511104$ degrees
of freedom at each time $t$.
This requires approximately $35.18$ Terabytes 
per temporal snapshot if we store the solution 
tensor $\bm f$ in a double precision 
IEEE 754 floating point format.
Several general-purpose algorithms have been developed 
to mitigate such an exponential growth of degrees of freedom, 
the computational cost, and the memory requirements. 
These algorithms include, e.g., sparse 
collocation methods \cite{Bungartz,Barthelmann,Akil} and 
techniques based on deep neural networks \cite{Raissi,Raissi1}.

In a parallel research effort that has its roots in quantum field 
theory and quantum entanglement, researchers have recently developed a new 
generation of algorithms based on tensor networks and low-rank 
tensor techniques to compute the solution of high-dimensional PDEs \cite{khoromskij,Bachmayr,parr_tensor,DektorVenturi2023,DektorVenturi2024,Heyrim2017}. A tensor network is a factorization of an entangled object such as a multivariate function or an operator into simpler objects, e.g., low-dimensional 
functions or separable operators, which are amenable to 
efficient representation and computation. A tensor network can 
be visualized in terms of graphs (see Figure \ref{fig:TN}). 
\begin{figure}[t]
\begin{center}
\includegraphics[scale=0.32]{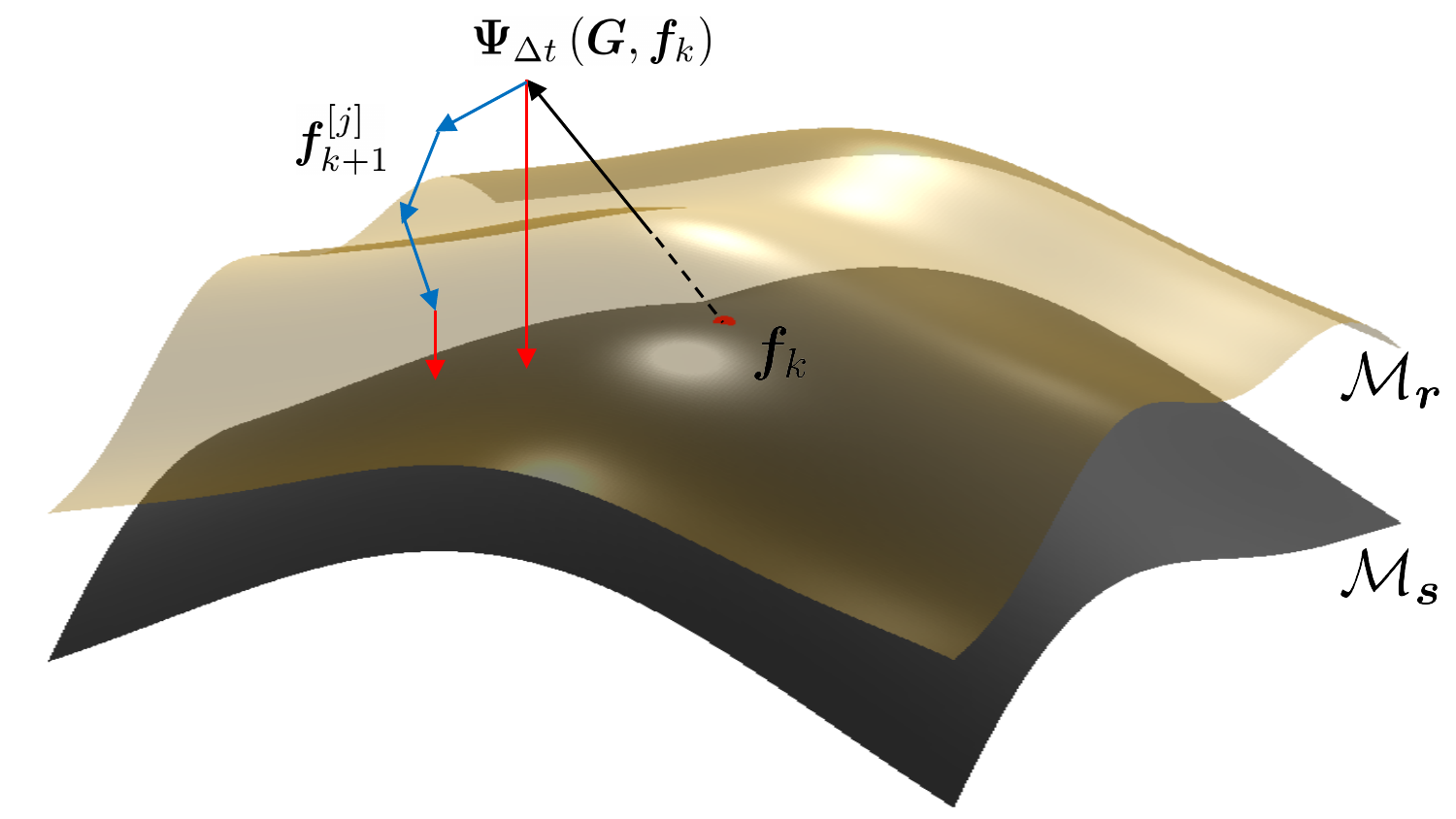}
\end{center}
\caption{
Sketch of implicit and explicit step-truncation 
integration methods. Given a tensor $\bm f_k$ with multilinear 
rank $\bm s$ on the tensor manifold ${\cal M}_{\bm s}$, we first 
perform an explicit time-step with any conventional 
time-stepping scheme. The explicit step-truncation integrator 
then projects $\bm \Psi_{\Delta t}(\bm G,\bm f_k)$ to a tensor manifold 
with rank $\bm s$. The multilinear rank $\bm s$ 
is chosen adaptively based on desired accuracy and 
stability constraints \cite{rodgers2020step-truncation}. 
On the other hand, the implicit step-truncation method 
takes $\bm \Psi_{\Delta t}(\bm G,\bm f_k)$ 
as input and generates a sequence of fixed-point 
iterates ${\bm f}^{[j]}$ shown as dots connected with 
blue lines. The last iterate is then projected 
onto a low rank tensor manifold, illustrated 
here also as a red line landing on
${\cal M}_{\bm s}$. This operation is 
equivalent to the {compression step} 
in the HT/TT-GMRES algorithm described 
in \cite{dolgov2013ttgmres}.}
\label{fig:intro-surf-compare}
\end{figure}
The vast majority of tensor algorithms currently available 
to approximate functions, operators and PDEs on tensor spaces 
is based on canonical polyadic (CP) decompositions \cite{Heyrim2017,parr_tensor,BoltzmannBGK2020,Beylkin2002}, 
Tucker tensors, or tensors corresponding to binary trees such as 
tensor train (TT) \cite{Dektor_dyn_approx,Bigoni_2016,OseledetsTT}
and hierarchical Tucker (HT)
tensors \cite{approx_rates,Grasedyck2018,Etter,h_tucker_geom}. 
A compelling reason for using binary tensor trees is that they allow 
to construct the tensor expansion by leveraging the spectral theory 
for linear operators, in particular the hierarchical 
Schmidt decomposition \cite{adaptive_rank,Dektor_dyn_approx,Griebel2019,rodgers2020stability,kato}.

\subsection{Step-truncation tensor methods}
\label{sec:step-truncation}
A new class of algorithms to integrate \eqref{eqn:ode} on 
a low-rank tensor manifold was recently proposed 
in \cite{rodgers2020step-truncation,rodgers2023implicit}.
These algorithms are known as {\em step-truncation methods}
and they are based on integrating the solution ${\bm f}(t)$ of the ODE \eqref{eqn:ode}
off the tensor manifold for a short time using any conventional explicit 
time-stepping scheme, and then mapping it back onto the manifold 
using a tensor truncation operation. 
To briefly describe these methods, let us 
discretize the ODE \eqref{eqn:ode} in time with a 
one-step method on an evenly-spaced temporal grid as 
\begin{align}
{\bm f}_{k+1} = 
{\bm \Psi}_{\Delta t}({\bm G}, {\bm f}_{k}),
\qquad {\bm f}_{0}={\bm f}(0),
\label{eqn:discrete-flow-intro}
\end{align}
where ${\bm f}_{k}$ denotes an approximation of 
${\bm f}(k\Delta t)$ for $k=0,1,\ldots$, 
and ${\bm \Psi}_{\Delta t}$ is an increment function.
To obtain a step-truncation integrator, we simply apply 
a truncation operator $\mathfrak{T}_{\bm s}(\cdot)$, i.e., 
a nonlinear projection onto a tensor manifold ${\cal M}_{\bm s}$
of HT (or TT) tensors with multilinear rank 
$\bm s$ \cite{uschmajew2013geometry},
to the scheme \eqref{eqn:discrete-flow-intro}  
(see Figure \ref{fig:intro-surf-compare}). 
This yields 
\begin{align}
\label{eqn:discrete-flow-trunc-intro}
{\bm f}_{k+1} = 
{\mathfrak T}_{\bm s}\left[
{\bm \Psi}_{\Delta t}({\bm G}, {\bm f}_{k})\right].
\end{align}
The need for tensor rank-reduction when iterating 
\eqref{eqn:discrete-flow-intro} can be easily understood 
by noting that tensor operations such as the application of 
an operator to a tensor and the addition between 
two tensors naturally increase tensor rank
\cite{kressner2014algorithm}.

We now present the rank-adaptive step-truncation 
algorithms used here to integrate FDEs 
on tensor manifolds. 
The first one is a truncated version of the Euler 
forward scheme, i.e., 
\begin{equation}
\label{adaptive-euler-method}
{\bm f}_{k+1} =
{\mathfrak{T}}_{{\bm r}}[{\bm f}_k +
\Delta t {\mathfrak{T}}_{{\bm s}}[{\bm G}( {\bm f}_k)]].
\end{equation}
We have shown that this scheme is convergent to first-order 
in $\Delta t$ provided the rank $\bm r$ is properly updated 
at each time step (see \cite{rodgers2020step-truncation} for details).
\begin{figure*}[t]
\centerline{
\includegraphics[scale=0.5]{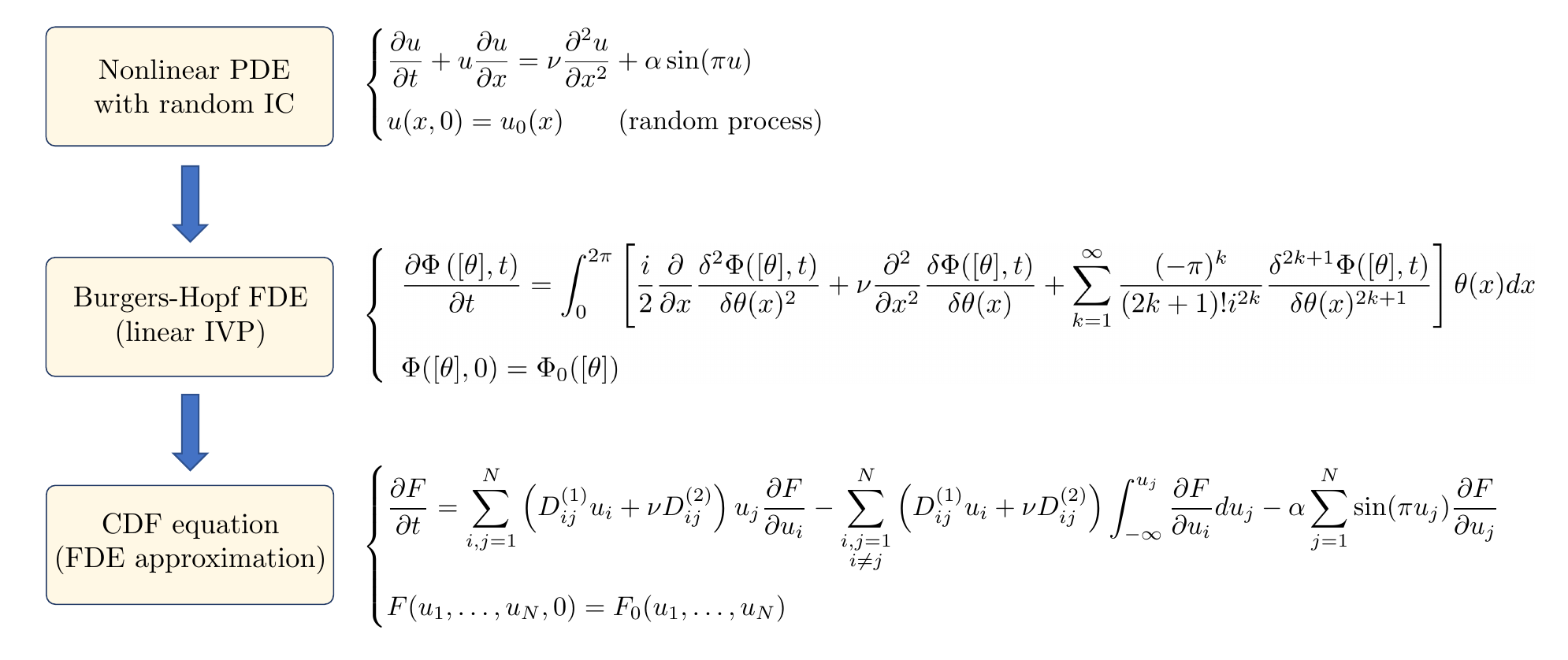}
}
\caption{Sketch of the logical flow to transform a random nonlinear initial value problem (IVP) for the Burgers equation into the (linear) Burgers-Hopf FDE, which is linear. The Burgers-Hopf FDE is subsequently approximated using the methods outlined in Section \ref{sec:FTE_to_hyperPDEs}. This yields an $N$-dimensional characteristic function equation of the form \eqref{hyper1}. By taking the inverse Fourier transform of such characteristic function equation and integrating it with respect to the phase variable we obtain the $N$-dimensional CDF equation shown above. The CDF equation is then solved using the proposed step-truncation tensor methods. }
\label{fig:fde_to_pde}
\end{figure*}
The second step-truncation scheme is derived from the explicit
midpoint rule, and it is defined as 
\begin{equation}
\label{adaptive-midpoint-method_2}
    {\bm f}_{k+1} =
    {\mathfrak{T}}_{{\bm \alpha}} \left [
    {\bm f}_{k} + \Delta t
    {\mathfrak{T}}_{{\bm \beta}}\left [
    {\bm G}\left (
    {\bm f}_{k} + \frac{\Delta t}{2}
    {\mathfrak{T}}_{{\bm \gamma}}(
    {\bm G}({\bm f}_{k}))
    \right )\right ] \right ], 
\end{equation}
where $\bm \alpha$, $\bm \beta$, and $\bm \gamma$ are once again selected adaptively as time integration advances. 

Both \eqref{adaptive-euler-method} and 
\eqref{adaptive-midpoint-method_2} are {\em explicit} 
step-truncation methods. These schemes are very simple 
to implement and have proven successful in integrating initial 
value problems for a variety of PDEs \cite{rodgers2020step-truncation}.
However, the combination of dimensionality, non-linearity and 
stiffness may introduce time-step restrictions which 
could make explicit time integration on tensor manifolds  
computationally infeasible.
To overcome this problem, we recently developed 
a new class of {\em implicit} rank-adaptive step-truncation 
algorithms for temporal integration PDEs on tensor manifolds \cite{rodgers2023implicit}.
These algorithms are based on an inexact Newton's 
method with tensor train GMRES iterations \cite{dolgov2013ttgmres} 
to solve the algebraic equation arising from time 
discretization on the tensor manifold.

To describe these methods, let us consider the standard Euler backward scheme  
\begin{equation}
\label{eqn:imp-euler}
{\bm f}_{k+1} = {\bm f}_{k} + \Delta t {\bm G}({\bm f}_{k+1}),
\end{equation}
and the associated root-finding problem 
\begin{equation}
\bm H_k(\bm f_{k+1})={\bm f}_{k+1} - {\bm f}_{k} -
\Delta t {\bm G}({\bm f}_{k+1})={\bm 0}.
\label{EB}
\end{equation}
If $\bm G$ is linear, e.g., for PDEs of the form \eqref{hyper1},
this reduces to a linear inversion problem
\begin{equation}
\left ({\bm I} - \Delta t{\bm G}\right )
{\bm f}_{k+1}={\bm f}_{k}.
\label{EB-linear}
\end{equation}
To solve the linear system \eqref{EB-linear} in a tensor format,
we apply the TT-GMRES method proposed by Dolgov in \cite{dolgov2013ttgmres}. 
For the implicit midpoint method, the logic
is identical. We apply the above steps to the
linear system
\begin{equation}
\left ({\bm I} - \frac{1}{2}\Delta t{\bm G}\right )
{\bm f}_{k+1}=
\left ({\bm I} + \frac{1}{2}\Delta t{\bm G}\right )\bm f_{k}.
\label{MP-linear}
\end{equation}
For the fully nonlinear case, one may couple the
TT-GMRES iteration to an inexact Newton method
and retain all convergence properties 
(see \cite[Appendix A]{rodgers2023implicit} for further details). 
In \cite{rodgers2023implicit} we have recently show that 
under mild conditions on the tensor truncation error, the implicit 
step-truncation Euler and implicit step-truncation midpoint 
converge, respectively, with order one and order two in $\Delta t$.

\section{Burgers-Hopf FDE}
\label{sec:applicationBurgers}

\begin{figure*}[t]
\centerline{\hspace{-3.cm}
\footnotesize 
Local Truncation Error 
}
\centerline{ 
\includegraphics[scale=0.65]{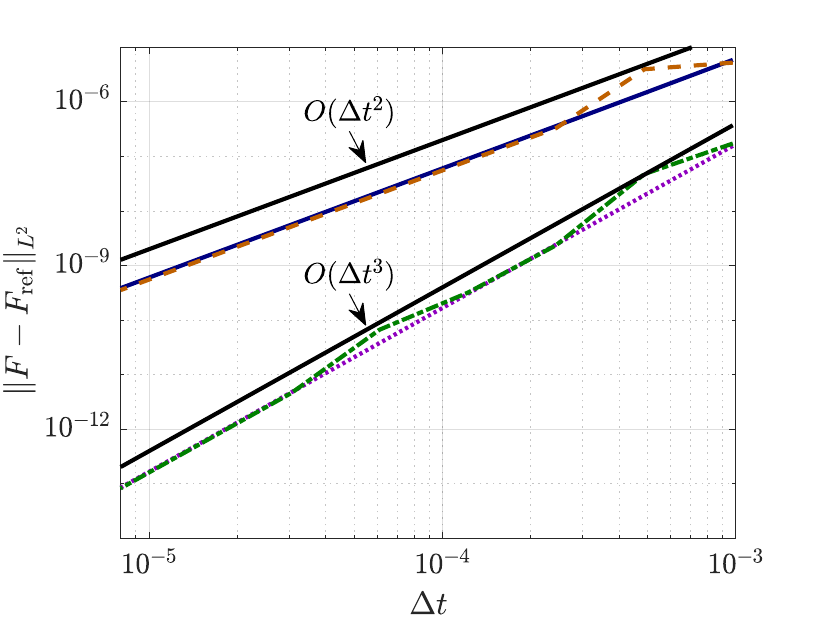}
 \includegraphics[scale=0.65]{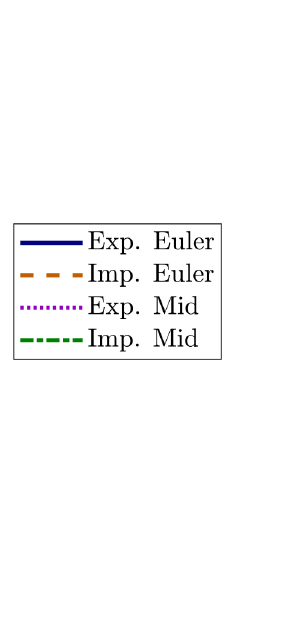}
}
\caption{
Local truncation errors (calculated as $O(\Delta t^{p+1})$ for
and order $p$ method)
of the proposed step-truncation tensor methods.
These errors are computed by comparing one time step of step-truncation method with its Richardson extrapolation. It is seen that  explicit/implicit Euler and 
explicit/implicit midpoint step-truncation methods have 
accuracy of order one and order two, respectively.
}
\label{fig:err-dim-8}
\end{figure*}

\begin{figure*}[ht]
\centering
\centerline{\hspace{2.0cm}
\footnotesize 
$t = 0.0$ \hspace{2.5cm}
$t = 0.04$ \hspace{2.5cm}
$t = 0.08$ \hspace{2.5cm}
$t = 0.12$ 
}
\centerline{\line(1,0){500}}
\vspace{0.1cm}
{\rotatebox{90}{\hspace{1.6cm}\rotatebox{-90}{
\hspace{0.1cm}
\footnotesize
MC CDF
\hspace{0.35cm}
}}}
\includegraphics[scale=0.47]{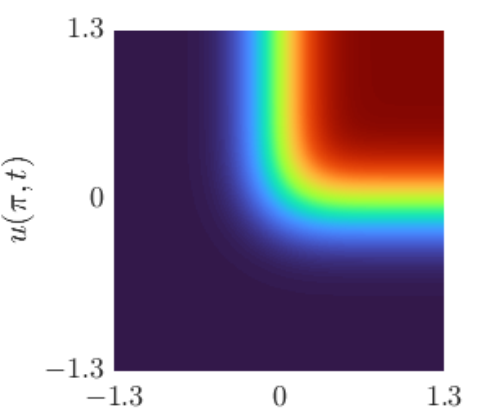}
\includegraphics[scale=0.47]{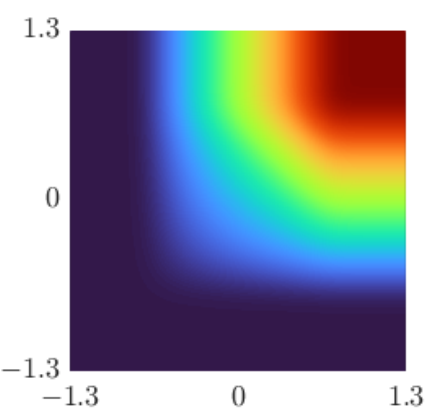}
\includegraphics[scale=0.47]{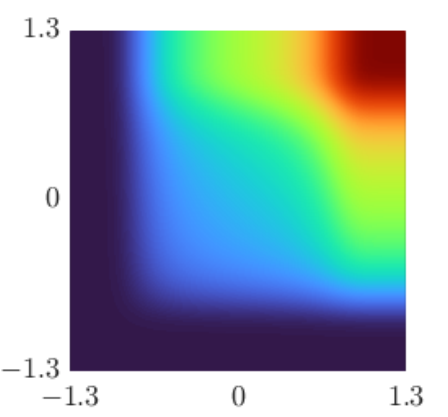}
\includegraphics[scale=0.47]{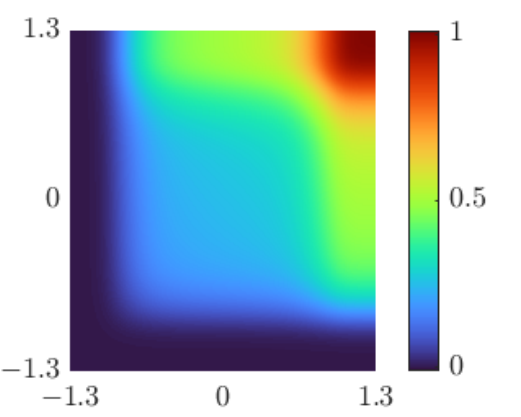}\\
{\rotatebox{90}{\hspace{1.6cm}\rotatebox{-90}{
\hspace{0.05cm}
\footnotesize
ST CDF
\hspace{0.5cm}
}}}
\includegraphics[scale=0.465]{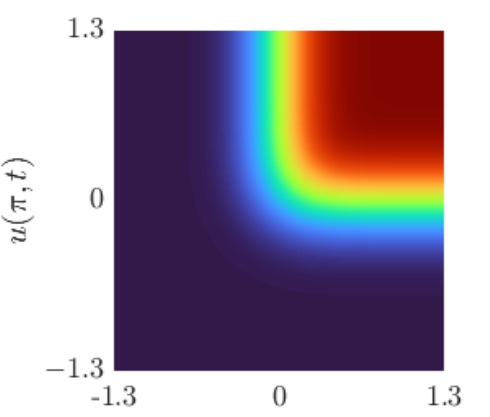}
\includegraphics[scale=0.47]{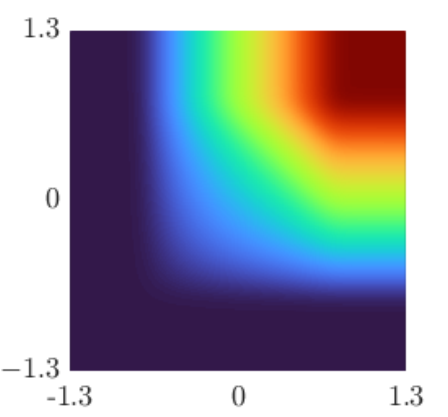}
\includegraphics[scale=0.47]{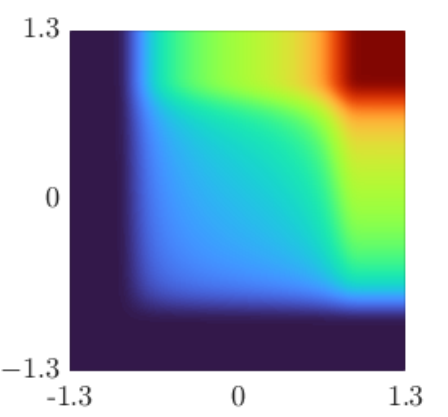}
\includegraphics[scale=0.47]{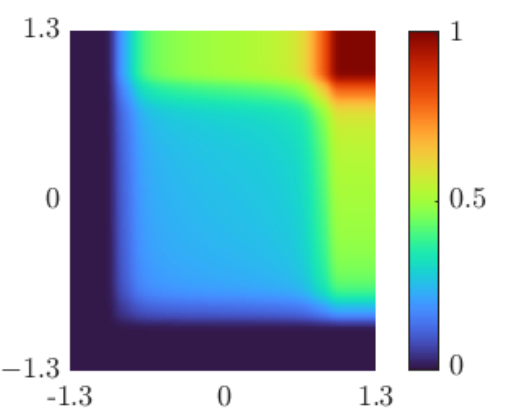}\\
{\rotatebox{90}{\hspace{1.6cm}\rotatebox{-90}{
\hspace{0.05cm}
\footnotesize
MC PDF
\hspace{0.35cm}
}}}
\includegraphics[scale=0.48]{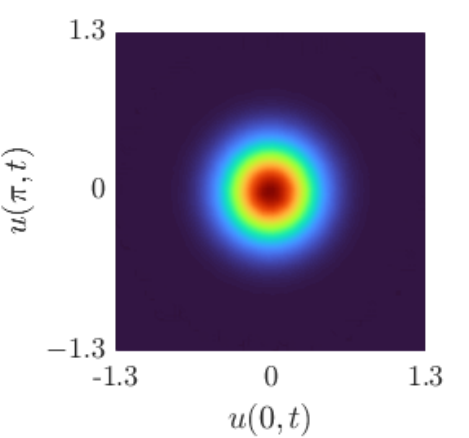}
\includegraphics[scale=0.47]{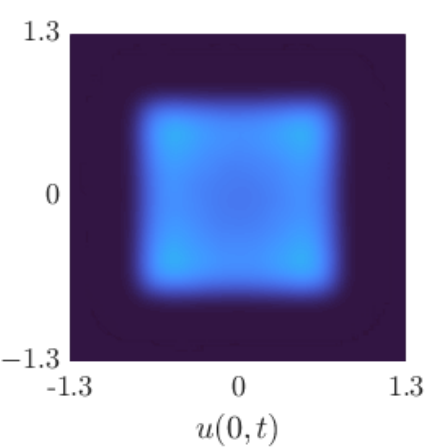}
\includegraphics[scale=0.47]{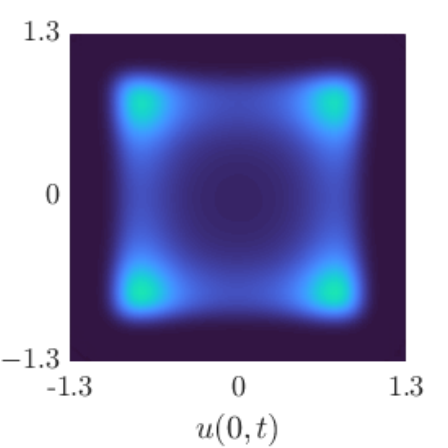}
\includegraphics[scale=0.47]{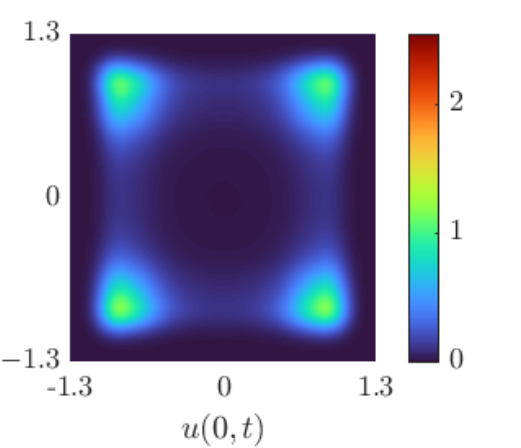}\\
\caption{
Two joint CDF of $u(0,t)$ and $u(\pi,t)$.
The CDF is computed by generating
numerical solutions to \eqref{eqn:cdf-burgers}
for $N=20$ using the proposed step-truncation tensor methods, 
and then marginalizing the solution in the remaining 18 variables.
We also show a Monte-Carlo estimate
of the joint CDF obtained by sampling $5\times 10^6$
solutions to \eqref{eqn:burgers-rxn}.}
\label{fig:time-snapshot-plot-marg-0-mid}
\end{figure*}

Consider Burger's equation on the unit circle $ \Omega=[0,2\pi]$
\begin{equation}
    \frac{\partial u}{\partial t} + u\frac{\partial u}{\partial x}
    = \gamma\frac{\partial^2 u}{\partial x^2} + R(u),
    \label{eqn:burgers-rxn}
\end{equation}
with a reaction term $R(u) = \alpha\sin(\pi u)$ and
random ($2\pi$-periodic) initial condition $u(x,0)$.
Let 
\begin{equation}
 \Phi([\theta(x)],t) =\mathbb{E}\left\{\exp\left(i\int_{0}^{2\pi}u(x,t)\theta(x)dx\right)\right\}   
\end{equation}
be the characteristic functional 
of the stochastic solution to \eqref{eqn:burgers-rxn}. It is straightforward to show that $\Phi([\theta(x)],t)$ satisfies the following Burgers-Hopf FDE
\begin{align}
\frac{\partial \Phi([\theta],t)}{\partial t}
=
&\int_{0}^{2\pi}
\Biggl [
\frac{i}{2}
\frac{\partial }{\partial x}
\frac{\delta^2 \Phi([\theta],t)}{\delta \theta(x)^2}
\nonumber
+\gamma \frac{\partial^2 }{\partial x^2}
\frac{\delta \Phi([\theta],t)}{\delta \theta(x)}\\
&+\alpha\sum_{k=1}^\infty \frac{(-\pi)^k}{(2k+1)!i^{2k}}
\frac{\delta^{2+1} \Phi([\theta],t)}{\delta \theta(x)^{2k+1}}
\Biggl ]\theta(x)dx,
\label{eqn:BHfde}
\end{align}
which admits a unique solution due to the uniqueness 
of the solution to the weak form of the Burgers equation 
(see Section \ref{sec:exist_and_unique}).

By approximating the FDE \eqref{eqn:BHfde} in terms of 
a high-dimensional PDEs of the form of \eqref{hyper1} and 
taking the inverse Fourier transform (in the sense of 
distributions) yields the Liouville-type linear 
hyperbolic conservation law
\begin{align}
\frac{\partial p}{\partial t} =
-\sum_{i=1}^{N}\frac{\partial}{\partial u_j}\Biggl[&
\sum_{j=1}^{N}
\Biggl (-u_iD_{ij}^{(1)}u_j+\nonumber \\
&\quad \delta_{ij}R(u_j)+
\gamma D_{ij}^{(2)}u_j
\Biggl )p
\Biggl ],
\label{eqn:pdf-burgers}
\end{align}
where $\delta_{ij}$ is the Kronecker delta and 
$p(u_1,u_2,\dots u_N,t)$ is the joint probability density function
of the solution to \eqref{eqn:burgers-rxn} at all evenly-spaced 
grid points $x_j= 2 \pi (j-1)/N$ ($j=1,\ldots, N$), 
i.e., $u_j=u(x_j,t)$.

In \eqref{eqn:pdf-burgers}, $D_{ij}^{(1)}$
and $D_{ij}^{(2)}$ are differentiation matrices. 
In our numerical experiments, we use the second-order 
finite difference matrix on a periodic domain.
From here, it is straightforward to write down a 
step-truncation tensor method to solve 
the high-dimensional PDE \eqref{eqn:pdf-burgers}.

However, this may not yield an accurate scheme 
due to the properties of the solution. In fact, it is well-known that 
PDF equations of the form \eqref{eqn:pdf-burgers} can develop shocks \cite{HeyrimPRS_2014}, or can have solutions converging to 
Dirac delta functions if the time asymptotic solution 
is deterministic. In these cases, it is convenient to first transform 
the kinetic equation into a cumulative distribution function 
(CDF) equation via an integral transform \cite{Daniele2013}, and then 
solve such CDF equation. This allows us to transform PDF shocks into CDF 
ramps, and PDF Dirac deltas into CDF shocks. In other words, 
the CDF equations produces solutions with features that 
are well-known to the numerical analysis community (shocks, ramps 
and solutions with range bounded in $[0,1]$) 
and therefore can be handled with numerical schemes for hyperbolic 
conservation laws \cite{Hesthaven2018}.

\begin{figure*}[t]
\centerline{ 
\includegraphics[scale=0.6]{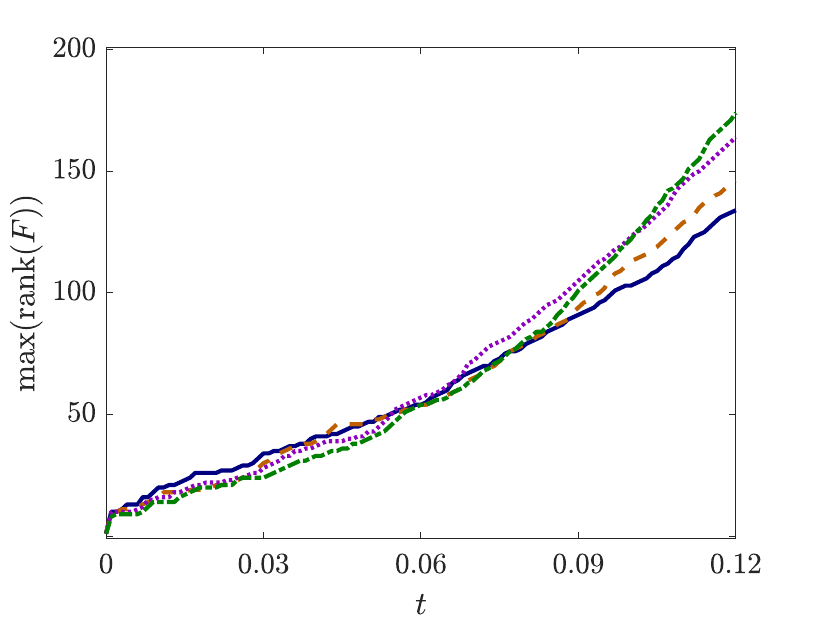}
\includegraphics[scale=0.6]{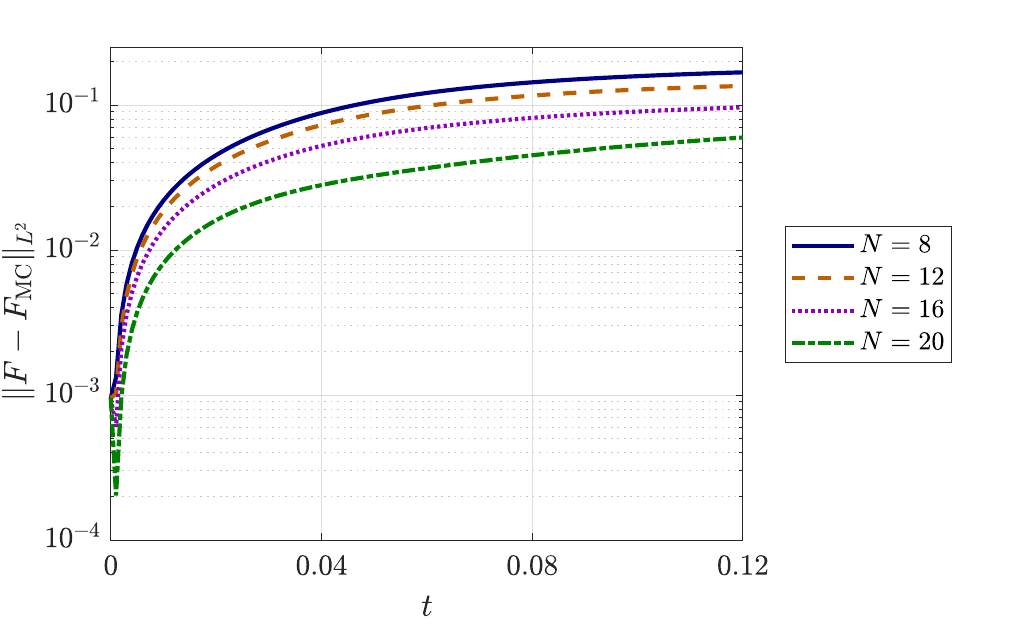}
}
\caption{
Left: Highest rank of the Tensor Train cores for 
tensor solutions of the CDF to equation \eqref{eqn:cdf-burgers}
computed with the explicit step-truncation midpoint method for 
dimension $N=\{8,12,16,20\}$.
Right: Convergence of the solution to the CDF equation 
\eqref{eqn:cdf-burgers} the Burgers-Hopf FDE as we increase 
the dimension $N$ from $8$ to $20$. The error plot represents 
represents the $L^2$ error between the joint CDFs shown in Figure
\ref{fig:time-snapshot-plot-marg-0-mid}.}
\label{fig:rank-err-highdim}
\end{figure*}

By integrating \eqref{eqn:pdf-burgers} with respect to the phase variables 
$\{u_1,\ldots,u_N\}$ and using we obtain the following CDF 
equation
\begin{align}
\nonumber
  \frac{\partial F}{\partial t}  &=
\sum_{i,j=1}^N
(u_iD_{ij}^{(1)}u_j -\gamma D_{ij}^{(2)}u_j-\delta_{ij}R(u_j))
\frac{\partial F}{\partial u_i}\\
&\quad -
\sum_{\substack{{i,j=1}\\
{i\neq j}}}^N
(D_{ij}^{(1)}u_i - \gamma D_{ij}^{(2)})
\int_{-\infty}^{u_j}
\frac{\partial F}{\partial u_i}d u_j\nonumber\\
&\quad -\alpha\sum_{j=1}^N \sin(\pi u_j)\frac{\partial F}{\partial u_j},
\label{eqn:cdf-burgers} 
\end{align}
where
\begin{equation}
    F(u_1,\ldots,u_n,t) =\int_{-\infty}^{u_1}\cdots \int_{-\infty}^{u_n} p(v_1,\ldots,v_n,t)dv_1\cdots dv_n.
\end{equation}
With the CDF $F(u_1,\ldots,u_N,t)$ available it is 
straightforward to compute the PDF
\begin{equation}
    p(u_1,\ldots,u_N,t)=\frac{\partial^n F(u_1,\ldots,u_N,t)}{\partial u_1\cdots\partial u_N}.
\end{equation}

To apply the proposed step-truncation tensor 
methods to equation \eqref{eqn:cdf-burgers}, 
we discretize the phase space in 
terms of a $N$-dimensional hypercube centered at the origin 
which covers the support of $p(u_1,u_2,\dots u_N,t)$. For the 
specific application at hand such hypercube is $\Sigma=[-1.3,1.3]^N$.
On the boundary of $\Sigma$, it is straightforward
to show that the derivative of $F$ normal to the hypercube
boundary is zero. Therefore we may apply outflow
boundary conditions, i.e., set $\partial F/\partial u_j=0$ 
at the boundary of $\Sigma$. 
We discretize each variable $u_i$ in $\Sigma$ on an 
evenly-spaced grid with $n_i=64$ points. Also, we approximate 
the integrals in \eqref{eqn:cdf-burgers} using the trapezoidal 
rule, and the partial derivatives with centered second-order 
finite differences. For the initial condition, 
we set the CDF to be a product of $N$ zero-mean 
Gaussian CDFs with standard deviation $\sigma=0.25$, i.e., 
\begin{equation}
F(u_1,\ldots,u_N,0)=\prod_{k=1}^N g(u_k),
\end{equation}
where 
\begin{equation}
    g(u)=\frac{1}{2} \left[1+\text{erf}\left(\frac{u}{\sigma\sqrt{2}}\right)\right]
\end{equation}
and $\text{erf}(x)$ is the standard error function.

\subsection{Numerical results}

To solve the high-dimensional CDF equation \eqref{eqn:cdf-burgers} 
approximating the Burgers-Hopf FDE \eqref{eqn:BHfde} 
we developed a high-performance (parallel) tensor code 
implementing the step truncation tensor methods discussed in 
Section \ref{sec:step-truncation}. 
The code can be found at the GitHub repository: 
\texttt{https://github.com/akrodger/paratt}.
        
We first perform a numerical study to verify 
temporal accuracy of order 1 and order 2 of the proposed 
step-truncation tensor methods applied to the 
high-dimensional CDF equation \eqref{eqn:cdf-burgers}
(see Figure \ref{fig:fde_to_pde}).
To stabilize explicit step-truncation methods, we add 
a numerical diffusion term proportional to the round-off 
error in $\Delta t$ to the right hand
side of \eqref{eqn:cdf-burgers}, in a similar manner to the
Lax-Wendroff method \cite{lax1958systems}.
To guarantee that we are accurately
capturing second and first order errors
in time, we must ensure that the semi-discrete
PDE we are solving for our reference solution
matches that of our step-truncation integrators.
To this end, we compute the local truncation errors
of each scheme by computing its difference
with a Richardson extrapolation
(see \cite[II.4]{HairerErnst1993SODE}) of it to 1 order
higher, i.e., we take two steps of size $\Delta t/2$
and combine it with a step of size $\Delta t$ to
create a more accurate estimate for the purposes
of comparing it to a single step of size  $\Delta t$.
To this end, we calculate the local truncation errors 
by computing the difference between the prediction of 
each scheme and a Richardson extrapolation of one order higher 
(as outlined in \cite[II.4]{HairerErnst1993SODE}). In other words, 
we take two steps of size $\Delta t/2$ and integrate them with 
a step of size $\Delta t$ to generate a more precise estimation 
for comparison with a single step of size $\Delta t$. The results 
are shown in Figure \ref{fig:err-dim-8}. 
It is seen that the explicit/implicit Euler and 
explicit/implicit midpoint step-truncation integrators have 
accuracy of order one and order two, respectively. 
To quantify the errors, we constructed a 
benchmark CDF via Monte-Carlo simulation. 
In practice, we sampled $5\times 10^6$ initial conditions $u(x,0)$
from the Gaussian distribution described above, 
computed the solution $u(x,t)$ of \eqref{eqn:burgers-rxn} using 
the Fourier pseudospectral method \cite{Hesthaven}, and 
then used a kernel density estimator (KDE) \cite{botev2010kernel}
to obtain the one-point and the two-points CDFs. 
This allows us to study convergence of the solutions
to the CDF equation \eqref{eqn:cdf-burgers} as we increase 
the dimension $N$ of the PDE. This convergence study is {\em crucial} 
for establishing convergence of the solution of the CDF equation 
\eqref{eqn:cdf-burgers} to the solution of the Burgers-Hopf FDE. 
For this analysis, we selected the explicit step-truncation 
midpoint tensor method as the temporal integrator, due to its
favorable performance and ease of implementation.

\begin{figure*}[t]
\centering
\centerline{
\hspace{0.5cm}
$t = 0.0$ \hspace{5.2cm}
$t = 0.04$
}
\includegraphics[scale=0.48]{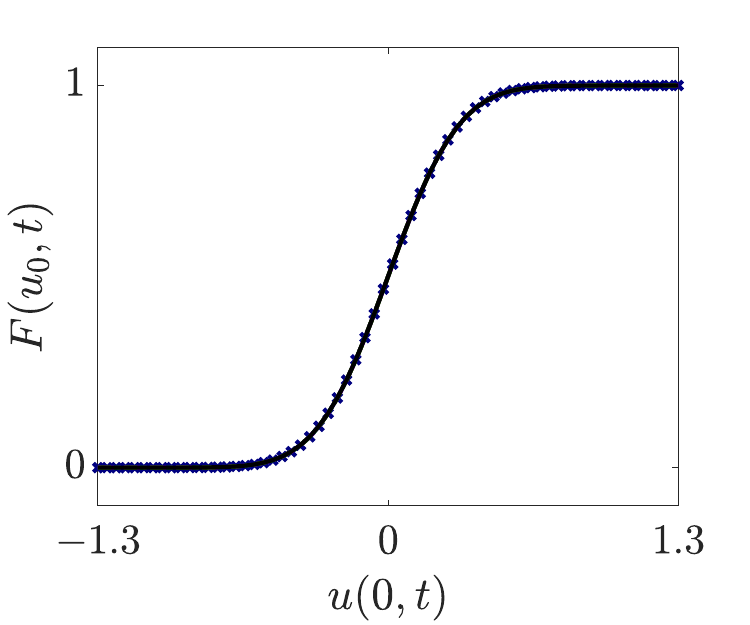}
\includegraphics[scale=0.48]{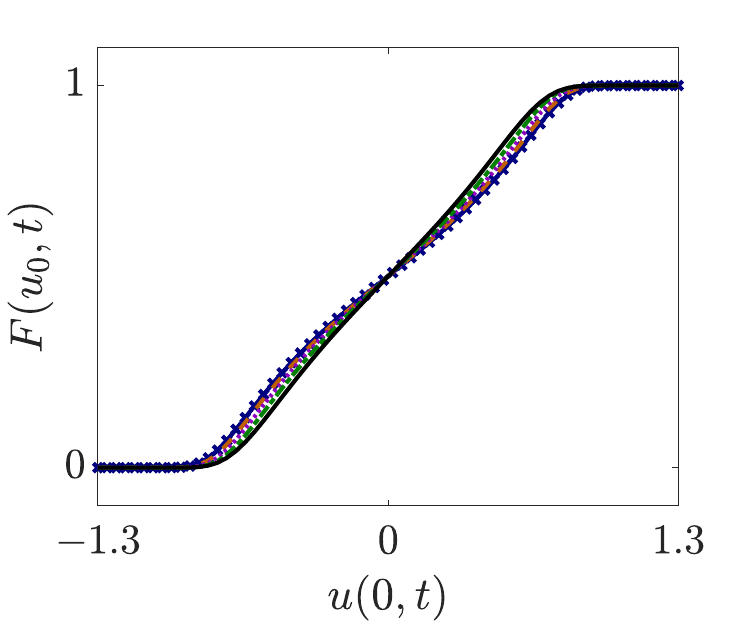}\\
\vspace{0.1cm}
\centerline{
\hspace{0.5cm}
$t = 0.08$ \hspace{5.2cm}
$t = 0.12$
}
\includegraphics[scale=0.48]{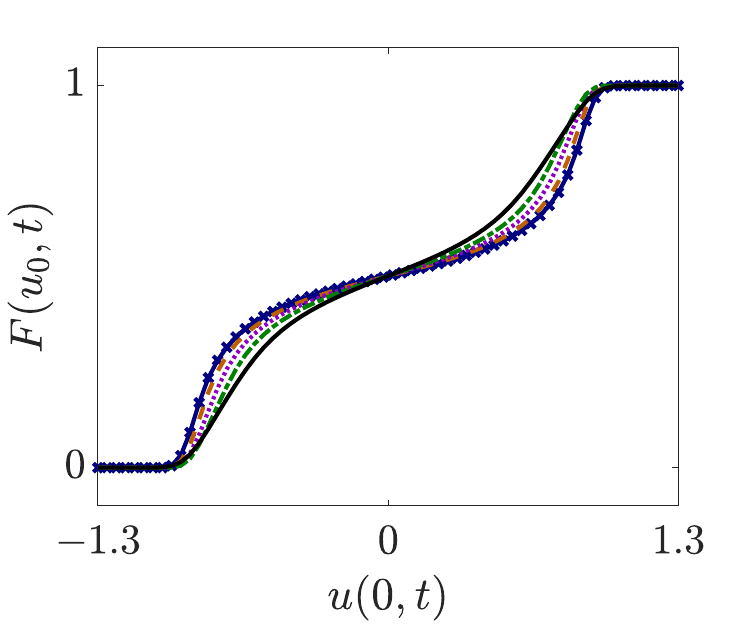}
\includegraphics[scale=0.48]{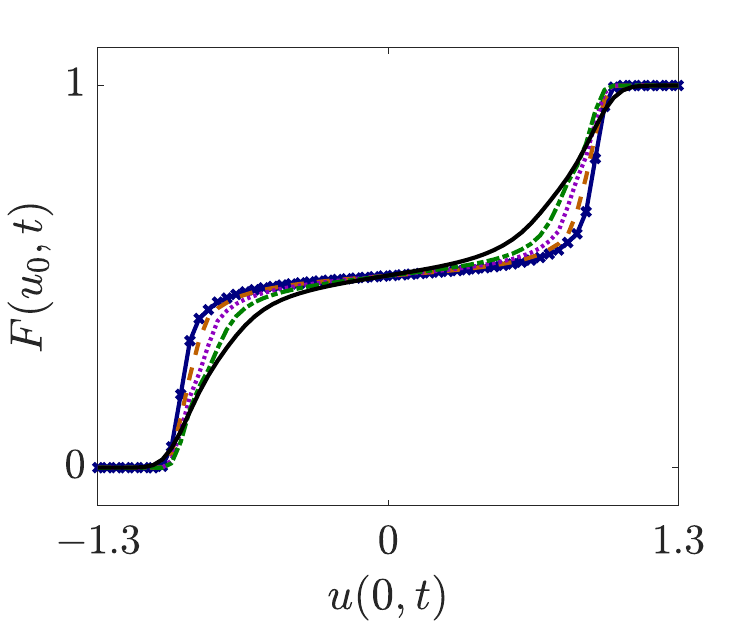}\\
\vspace{0.1cm}
\centerline{ 
\includegraphics[scale=0.48]{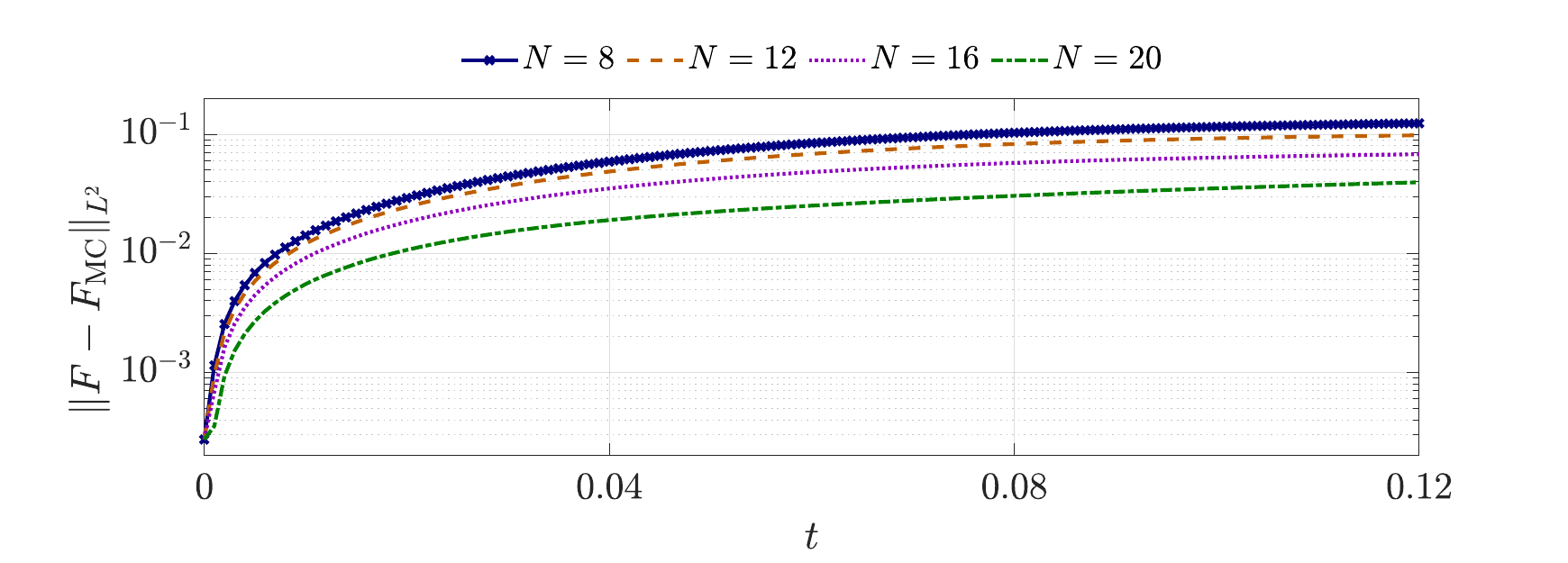}
\hspace{0.1cm}
}
\caption{CDF of $u(0,t)$ computed by solving 
equation \eqref{eqn:cdf-burgers} with the proposed 
step-truncation tensor method.  
It is seen that as we increase the dimension $N$ 
(number of independent variables in the CDF equation 
\eqref{eqn:cdf-burgers}) the tensor solution 
converges to the Monte-Carlo benchmark (black continuous line).}
\label{fig:time-snapshot-plot-marg-0}
\end{figure*}

In Figure \ref{fig:time-snapshot-plot-marg-0-mid}
plot the joint CDF of $u(0,t)$ and $u(\pi,t)$ at 
$t=\{0,0.04,0.08,0.012\}$ computed with Monte Carlo, 
and the proposed step-truncation midpoint tensor 
method applied to a 20-dimensional ($N=20$) CDF equation 
\eqref{eqn:cdf-burgers}. Clearly there is visual agreement 
between the results form the tensor simulation and the results form 
Monte-Carlo simulation. In Figure \ref{fig:rank-err-highdim} 
we demonstrate convergence of the solution to the CDF equation 
\eqref{eqn:cdf-burgers} the Burgers-Hopf FDE as we increase 
the dimension $N$ from $8$ to $20$. The error plot in 
Figure \ref{fig:rank-err-highdim} represents the $L^2$ 
error between the joint CDFs shown in Figure
\ref{fig:time-snapshot-plot-marg-0-mid}.
Runtime increases considerably with dimension though, 
especially for moderately high solution ranks. 
For example once rank surpasses $100$, the $N=20$ 
case requires approximately 45 minutes per time step 
on an Intel i9-7980xe workstation. The storage of the 
CDF tensor solution in 20 dimensions requires approximately 
690 megabytes per temporal snapshot for tensor ranks 
of the order of $150$. This is a significant compression 
of a multidimensional double precision floating point 
array of size $64^{20}$, which would normally require 
approximately $10^{22}$ petabytes. That is to say that we capture
$O(10^{-2})$ accuracy of the data while achieving a
data compression ratio of $1.44\times10^{26}\%$ and
data space savings of approximately $100\%$. Of course,
690 megabytes is not negligible on practical computing
systems. It is essentially nothing compared to the
decompressed data though.
In Figure \ref{fig:time-snapshot-plot-marg-0}
we compare the CDF of $u(0,t)$ computed by solving 
the CDF equation \eqref{eqn:cdf-burgers} with 
$N=\{8,12,16,20\}$ and then marginalizing, with a 
benchmark CDF obtained by Monte-Carlo simulation of \eqref{eqn:burgers-rxn}. 
It is seen that as we increase the dimension $N$ the 
marginal of the CDF tensor solution approaches 
the Monte-Carlo benchmark.
While the proposed numerical tensor schemes are accurate and provably 
convergent, they may not preserve important properties of the 
CDF solution, such as positivity, monotonicity, and range in $[0,1]$. 
The development of a structure-preserving tensor integration 
schemes is indeed a rather unexplored research area.
One possibility is to enforce structure at the level 
of the tensor truncation operation, i.e., perform 
rounding conditional to some required property. For instance, if
we are interested in enforcing non-negativity in 
a step-truncation tensor scheme, then we can simply replace 
the truncation operator based 
on recursive QR decomposition (see Appendix \ref{app:A}), with 
the distributed non-negative tensor truncation discussed 
in \cite{Bhattarai}. Another approach to enforce structure 
can be built directly at the level of the tensor equations. 
This approach was recently proposed by Einkemmer {\em et al.} 
in \cite{Einkemmer1,Einkemmer2,Einkemmer3} in 
the context of Vlasov and Vlasov-Poisson kinetic equations. 
The resulting dynamical low-rank algorithm conserves mass, 
momentum, and energy.

\section{Summary}
\label{sec:conclusion}

The numerical simulation of functional differential equations (FDEs) remains an open problem in many respects.  
In this paper, we developed approximation theory and new high-performance computational algorithms to solve FDEs on tensor manifolds. Our approach involves initially approximating the given FDE with a high-dimensional partial differential equation (PDE) using the functional approximation methods described in \cite{VenturiSpectral, venturi2018numerical, Ganbo2021}. Subsequently, we compute the numerical solution of such PDE using parallel rank-adaptive step-truncation methods \cite{rodgers2020step-truncation, rodgers2023implicit}. This is a novel and effective strategy for computing numerical solutions to FDEs, which paves the way for further research and practical applications across various branches of physics and engineering. We demonstrated the effectiveness of the proposed new approach through its application to the Burgers-Hopf FDE, which governs the characteristic functional of the stochastic solution to the Burgers equation evolving from a random initial state. 
Our numerical results and convergence studies demonstrate that, despite their computational complexity, tensor methods allow us to compute accurate 
solutions to the Burgers-Hopf FDE.

\begin{acknowledgments}
This research was supported by the U.S. Army Research Office grant W911NF1810309. Dr. Rodgers would like to acknowledge the Transformational Tools and Technologies $(\text{T}^3)$ Project for partially funding his Ph.D. research through the NASA Pathways internship program.
\end{acknowledgments}

\appendix
\section{Tensor algorithms}
\label{app:A}

In this Appendix, we describe the high-performance (parallel)
tensor algorithms we developed to solve high-dimensional PDEs on 
tensor manifolds, in particular the CDF equation discussed 
in Section \ref{sec:applicationBurgers}. 
To this end, we begin with a brief review of the algebraic 
representation of the tensor train format \cite{Oseledets2010}.

\subsection{Tensor train decomposition}
A tensor in ${\bm f}\in{\mathbb R}^{n_1\times n_2 \times \cdots \times n_d}$
is in a tensor train (TT) format if there is an array 
of positive integers (called the TT rank)
${\bm r}=(r_0,r_1,r_2,\dots ,r_{d-1},r_d)$ with
$r_0=1=r_d$ and a list of order 3 tensors
(called the TT cores)
${\bm C}=({\bm C}_1,{\bm C}_2,\dots,{\bm C}_d)$
where ${\bm C}_k\in{\mathbb R}^{r_{k-1}\times n_k \times r_{k}}$
such that the entries of ${\bm f}$ may be written
as the iterated matrix product
\begin{align}
{\bm f}[i_1,i_2,\dots,i_d]
=
{\bm C}_{1}[1,&i_1,:]
{\bm C}_{2}[:,i_2,:]\cdots \nonumber \\
&{\bm C}_{d-1}[:,i_{d-1},:]
{\bm C}_{d}[:,i_{d},1],
\end{align}
where the colon operator indicates multidimensional array slicing, 
akin to its usage in MATLAB.
Though rather involved at first glance, the above expression
may be derived by writing down a multivariate function
series expansion by the method of separation of variables,
rearranging the expression into a finite sequence
of infinite matrix
products, truncating the series so that the matrix
products are finite, then discretizing in space so that
the multivariate function is sampled on a tensor product grid.
From this perspective, it becomes apparent that the
tensor train cores represent a two- dimensional array of functions
of a single variable and the ranks are the number of functions
present in a series expansion approximation. A more
concrete exploration of this perspective is given in
\cite{dektor2020dynamically,Alec2020}.

Due the the iterated matrix product definition of the format,
it becomes apparent that a number of arithmetic operations
may be represented in the format by producing a new TT
tensor with different ranks. In particular if
${\bm f},{\bm g}$ have TT core lists ${\bm C},{\bm D}$
then by a simple inductive argument, we have
\begin{align}
\nonumber
{\bm f}&[i_1,i_2,\dots,i_d]+{\bm g}[i_1,i_2,\dots,i_d]=\\
\nonumber
&\quad\begin{bmatrix}
{\bm C}_{1}[1,i_1,:] \ \big | \   {\bm D}_{1}[1,i_1,:]
\end{bmatrix}
\begin{bmatrix}
{\bm C}_{2}[:,i_2,:] & \big |  & {\bm 0} \\
{\bm 0}              & \big |  & {\bm D}_{2}[:,i_2,:]
\end{bmatrix}
\\
&\cdots
\begin{bmatrix}
{\bm C}_{d-1}[:,i_{d-1},:] & \big |  & {\bm 0} \\
{\bm 0}              & \big |  & {\bm D}_{d-1}[:,i_{d-1},:]
\end{bmatrix}
\begin{bmatrix}
{\bm C}_{1}[:,i_d,1] \\   {\bm D}_{1}[:,i_d,1]
\end{bmatrix}.
\end{align}

\begin{algorithm}[t]
\caption{Left Orthogonalization}\label{alg:tt-l-orthog}
\KwData{A tensor $\bm f$ in TT format with cores
$({\bm C}_1,{\bm C}_2,\dots,{\bm C}_d)$.}
\KwResult{A tensor ${\bm g}={\bm f}$ in TT format with
left orthogonal cores.}
\vspace{5mm}
Reserve memory for each core
${\bm D}_c$ with sizes of ${\bm C}_c$.\\
${\bm D}_1 = {\bm C}_1$\\
\For{$c=1,2,\dots ,d-1$}{
	$[{\bm Q}_c,{\bm R}_c] = QR({\mathtt V}({\bm D}_c))$\\
	${\bm H}_{c+1} = {\bm R}_c{\mathtt H}({\bm C}_{c+1})$\\
	${\bm D}_{c}[:] = {\bm Q}_c[:]$\\  
	${\bm D}_{c+1}[:] = {\bm H}_{c+1}[:]$\\  
}
Set the cores of $\bm g$ as 
$({\bm D}_1,{\bm D}_2,\dots,{\bm D}_d)$
\end{algorithm}

Thus the sum ${\bm w} = {\bm f}+{\bm g}$ is a TT tensor with
cores obtained by concatenating those of ${\bm f}$ and ${\bm g}$.
Scalar multiplication is straightforward,
just scale any of the cores of the tensor.
\begin{align}
\alpha
{\bm f}[i_1,i_2,\dots,i_d]
=
(\alpha
{\bm C}_{1}&[1,i_1,:] )
{\bm C}_{2}[:,i_2,:]
\cdots\nonumber\\
&{\bm C}_{d-1}[:,i_{d-1},:]
{\bm C}_{d}[:,i_{d},1].
\end{align}
In order to perform more sophisticated
operations such as tensor truncations, we must frequently
reshape the tensor cores into lower dimensional arrays.
Two particularly useful reshapings are the horizontal flattening
\begin{align*}
{\mathtt H}:{\mathbb R}^{a\times b\times c}
&\longrightarrow
{\mathbb R}^{a\times(bc)}\\
{\bm C}[i,j,k] &\longmapsto {\bm H}[i,j+bk],
\end{align*}
and the vertical flattening,
\begin{align*}
{\mathtt V}:{\mathbb R}^{a\times b\times c}
&\longrightarrow
{\mathbb R}^{(ab)\times c}\\
{\bm C}[i,j,k] &\longmapsto {\bm V}[i+aj,k].
\end{align*}
The above maps are two dimensional array analogs of
the vectorization of a tensor, in which we have the
coordinates are listed as
${\bm C}[i,j,k] = {\bm v}[i+aj+abk]$. Note that
when the entries of a tensor are stored contiguously
in computer memory, no memory movement or copying
needs to be done to interpret a tensor as its flattening
or vectorization. To undo the flattening
in an algorithm, we denote copying
all the entries as  ${\bm C}[:] = {\bm H}[:]$
or ${\bm C}[:] = {\bm V}[:]$.

\begin{algorithm}[t]
\caption{Truncation of a TT tensor}\label{alg:tt-r-trunc}
\KwData{A tensor $\bm f$ in TT format with cores
$({\bm C}_1,{\bm C}_2,\dots,{\bm C}_d)$
and a desired accuracy $\varepsilon$.}
\KwResult{A tensor ${\bm g}$ so that
$\|{\bm g}-{\bm f}\|\leq \varepsilon$
in TT format with
right orthogonal cores.}
\vspace{5mm}
Reserve memory for each core
${\bm D}_c$ with sizes of ${\bm C}_c$.\\
Set $({\hat {\bm C}}_1,{\hat {\bm C}}_2,
\dots,{\hat {\bm C}}_d)$ as orthogonal by Alg.
\ref{alg:tt-l-orthog}.\\
${\hat \varepsilon} = \varepsilon/\sqrt{d-1}$\\
${\bm D}_1 = {\hat {\bm C}}_1$\\
\For{$c=d,d-1,\dots ,2$}{
	$[{\bm U}_c,{\bm \Sigma}_c,{\bm V}_c^\top] =
			SVD({\mathtt H}({\bm D}_c),{\hat \varepsilon})$
			\ \  (Truncated SVD)\\
	${\bm H}_{c-1} = {\mathtt V}({\hat {\bm C}}_{c-1})
				{\bm U}_c {\bm \Sigma}_c $\\
	${\bm D}_{c}[:] = {\bm U}_c[:]$\\  
	${\bm D}_{c-1}[:] = {\bm H}_{c-1}[:]$\\  
}
Set the cores of $\bm g$ as 
$({\bm D}_1,{\bm D}_2,\dots,{\bm D}_d)$
\end{algorithm}
The general procedure of truncating the ranks of a TT
tensor comes in two phases. The first phase is to
collect all the norm of the tensor into a single
core through a sequence of matrix factorizations, leaving
all other cores to have Frobenius norm 1. Then a sequence
of SVDs are applied, multiplying a non-unitary factor
from a tensor train core to its neighboring cores.
The first phase is called orthogonalization of
a TT tensor and we now present an algorithm for it.

There are two variants of orthogonalization we will
discuss. One is Left-to-Right orthogonalization and
the other Right-to-Left orthogonalization.
They are the same algorithm mirrored across the
middle core of the tensor train and each step in one
is a transpose of its counterpart in the other.
A TT tensor is left orthogonal if
for $c=1,2,\dots, d-1$, we have
${\mathtt V}({\bm C}_c)^\top {\mathtt V}({\bm C}_c) =
{\bm I}$. Similarly, a TT tensor is right orthogonal
if for $c=2,3,\dots, d$, we have
${\mathtt H}({\bm C}_c){\mathtt H}({\bm C}_c)^\top  =
{\bm I}$. One may transform the cores of a TT tensor
to be left orthogonal while the full tensor is unchanged
by applying iterated matrix factorizations.
All we must do is sequentially flatten, orthogonalize
by QR decompoosition,
then unflatten each core,  until the final one is
reached.
This process is summarized in Algorithm
\ref{alg:tt-l-orthog}.
For the right orthogonal case, the procedure
is largely identical, though we use
LQ factorization, the lower-triangular
orthogonalizing decomposition,  instead.

We now describe one algorithm for truncation of a tensor
train given an orthogonalized TT tensor. The process is
essentially the same as orthogonalization, we apply
a sequence of orthogonal matrix factorizations to each core,
multiplying one of the terms into that core's left or
right neighbor. We then replace the current core with
the entries of an orthogonal factor. The difference in
this algorithm is that we apply an SVD, rather than QR
or LQ factorization. The term ``truncation'' comes from truncating
the singular value decomposition
series expansion of a matrix up to a desired
tolerance. This truncation is described
as the rounding algorithm in \cite{oseledets2011tensor}.
We present the right orthogonal variant of it here as Algorithm
\ref{alg:tt-r-trunc}.
There is no clear general rule for preferring one
left or right variants truncation algorithm. They have identical
computational costs and similar formulations. It is possible
that differing runtime will occur based on row-major versus
column major layout of the matrices, though this is
highly dependent
on the specific memory mapping and workspaces used in the SVD,
QR, and LQ calls.

\subsection{Parallel algorithms for tensor train arithmetic}
Among the arithmetic operations performed on the tensor
train during the evolution of the solution to a differential
equation, the most costly is the truncation. This is due to the the fact
that it requires data access and editing of every core twice.
Due to the sequential nature of the loops in Algorithms
\ref{alg:tt-l-orthog} and \ref{alg:tt-r-trunc} we see that parallelizing
the inner loop would be most effective.
To this end, we follow  the approach of \cite{daas2020parallel}.
Our extension of the algorithms presented in \cite{daas2020parallel} 
can be found at the GitHub repository: 
\texttt{https://github.com/akrodger/paratt}.
Since the rank of the tensor train is expected to change frequently during program
execution, we must split the memory stored
on a total of $P$ compute nodes in a manner independent
of the tensor rank.

To this end, we introduce a ${P\times d}$
memory partition matrix ${\bm M}$
with positive integer values so that for each core 
${\bm C}_k\in
\mathbb{R}^{r_{k-1}\times n_k \times r_k}$, we have
$\sum_{p=1}^P {\bm M}[p,k] = n_k$.
This matrix describes a block-tensor storage layout
for the tensor train core list. More precisely,
for each core ${\bm C}_k$, we define an array of
core blocks $({\bm C}_k^{1},{\bm C}_k^{2},\dots ,{\bm C}_k^{P})$
with array sizes defined
via ${\bm C}_k^{p}\in\mathbb{R}^{r_{k-1}\times
{\bm M}[p,k]  \times r_k} $ so that
\begin{equation}
{\bm C}_k[i,:,j] = 
\begin{bmatrix}
\label{eqn:par-split}
{\bm C}_k^{1}[i,:,j]\\
{\bm C}_k^{2}[i,:,j]\\
\vdots \\
{\bm C}_k^{P}[i,:,j]
\end{bmatrix}.
\end{equation}
Each distributed memory compute node with
index $p$ stores the core list
$({\bm C}_1^{p},$
${\bm C}_2^{p},$
$\dots ,$
${\bm C}_d^{p})$.
This memory layout is designed so that
\begin{align*}
{\mathtt V}({\bm C}_k) &= 
\begin{bmatrix}
{\mathtt V}({\bm C}_k^{1})\\
{\mathtt V}({\bm C}_k^{2})\\
\vdots \\
{\mathtt V}({\bm C}_k^{P})
\end{bmatrix}
\quad
\text{and}\\
{\mathtt H}({\bm C}_k) &= 
\begin{bmatrix}
{\mathtt H}({\bm C}_k^{1})
& | &
{\mathtt H}({\bm C}_k^{2})
& | &
\cdots 
& | &
{\mathtt H}({\bm C}_k^{P})
\end{bmatrix}.
\end{align*}
We may therefore compute a flattening of the cores
without any memory movements or cross-node communications
by reinterpreting the array storage offsets in column
major layout. Additionally, sums and scalar multiplications
of tensors may also be computed in parallel without communications.

This parallel layout also allows for straightforward
applications of finite difference stencils. Consider a finite
difference stencil which requires $s$ many points where the number of
1 dimensional ghost cells required is $g=(s-1)/2$.
Let ${\bm D}_{FD}$ be the linear operator for this stencil.
To apply a partial
derivative in variable $k$ to a TT tensor, we replace
the entries of core ${\bm C}_k$ with the entries of
${\bm D}_{FD}$ applied to equation
\eqref{eqn:par-split} for each $i,j$. To perform the
parallel version of this stencil, we first perform a
nonblocking communication to deliver the ghost cells of
compute node $p$ to their neighboring
nodes $p\pm 1$. For nodes $p=1$ or $p=P$, we instead rely on 
discrete boundary condition formulas, communicating again
between nodes $p=1$ and $p=P$ if periodic boundary conditions
are required. After boundary condition communications,
every processor has an appropriate local version of
the finite difference stencil, which is identical save
for the definition of the ghost cells.

Numerical integration is similar, though instead of sharing
ghost cells of compute node $p$ with compute node $p\pm 1$,
we need only pass the data with compute node $p+1$, due to
the one-sided stencil of a cumulative summation formula.

\subsection{Parallel tall-skinny QR factorization}
In order to truncate a TT tensor with the described
block memory layout, we must perform a distributed memory
QR factorization on a matrix of the form
\begin{equation*}
{\bm A} = 
\begin{bmatrix}
{\bm A}_{1}\\
{\bm A}_{2}\\
\vdots \\
{\bm A}_{P}
\end{bmatrix},
\end{equation*}
or equivalently, perform an LQ factorization on the
horizontally concatenated transpose of $\bm A$.
For mathematical simplicity, we only present the 
QR variant, though our GitHub code available at \texttt{https://github.com/akrodger/paratt}
has both variants.
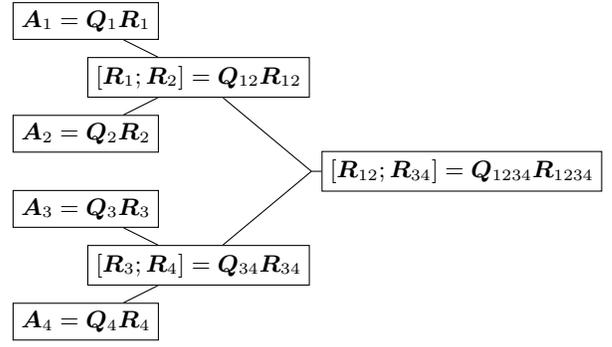
\begin{figure}
\begin{center}
\begin{tikzpicture}[every node/.style={rectangle,draw,inner sep=3pt}]
\node {{$[{\bm R}_{12};{\bm R}_{34}] = {\bm Q}_{1234}{\bm R}_{1234}$}} 
child[level distance=20mm,grow=left] {
child[level distance=15mm,sibling distance=25mm] {
		node{{$[{\bm R}_1;{\bm R}_2] = {\bm Q}_{12}{\bm R}_{12}$}} 
		    child[sibling distance=15mm]  {
		             node{${\bm A}_1 = {\bm Q}_1{\bm R}_1$}
		          }
		    child[sibling distance=15mm]  {
		             node{${\bm A}_2 = {\bm Q}_2{\bm R}_2$}
		          }
      } 
child[level distance=15mm,sibling distance=25mm] {
		node{{$[{\bm R}_3;{\bm R}_4] = {\bm Q}_{34}{\bm R}_{34}$}} 
		    child[sibling distance=15mm]  {
		             node{${\bm A}_3 = {\bm Q}_3{\bm R}_3$}
		          }
		    child[sibling distance=15mm]  {
		             node{${\bm A}_4 = {\bm Q}_4{\bm R}_4$}
		          }
      } 
};
\end{tikzpicture}
\end{center}
\caption{
Parallel TSQR algorithm with 4 compute notes.
A binary tree is formed which outlines the communication
pattern for calculating the ${\bm R}$ factor. At each level, the
node with larger ID sends its ${\bm R}$ factor to the
tree sibling it connects to. The sibling then does a QR
factorization with a concatenated pair of child ${\bm R}$ factors.
This process repeats until the final ${\bm R}$ factor is found at
the root of the tree stored on
processor $p=1$. We then broadcast this matrix
to all other processors. To get the ${\bm Q}$
factor, we store the orthogonal ${\bm Q}_i$
for each tree node $i$ and multiply
the parent's ${\bm Q}_i$ on the right side
of the child's ${\bm Q}$ factor,  traversing the tree
from the root to the leaves. This
results in a distributed memory
orthogonal matrix
${\bm Q} =[{\bm Q}_1;{\bm Q}_2;{\bm Q}_3;{\bm Q}_4].$
}
\label{fig:tsqr-tree}
\end{figure}
One such factorization amenable to this
layout is the Tall-Skinny QR (TSQR) factorization
\cite{demmel2012communication}.
This algorithm uses a binary tree structure to
define a communication pattern for decomposing the
QR factorization of $\bm A$ into a collection of
smaller QR factorizations, avoiding communications if
possible.
By factoring the blocks into their own QR factorizations,
we see that the Q factor may be presented as a sequence of
products of block diagonal matrices. This process
is then applied to a recursive binary tree to track
the nesting of the orthogonal Q factor.
The final result is an orthogonal matrix
${\bm Q} = [{\bm Q}_1;{\bm Q}_2;\dots ;{\bm Q}_P]$
distributed in the same
block layout as $\bm A$
and an upper triangular matrix $\bm R$ copied 
across all compute nodes.
See Figure \ref{fig:tsqr-tree} for
a visual representation of the binary tree storage structure
for the case $P=4$.
The transposed version of this algorithm is called the
Wide-Fat LQ (WFLQ) factorization and has the same
tree data storage structure, but every matrix factorization
is transposed.

\subsection{Parallel orthogonalization and truncation}

We now introduce a parallelization of Algorithms
\ref{alg:tt-l-orthog} and \ref{alg:tt-r-trunc}.
It can be seen that parallelizing the orthogonalization
algorithms is as simple as replacing the QR factorizations
with their TSQR variants. However, the truncation is
requires a bit of manipulation. In order to parallelize this,
we first note the following relationship
of the SVD and QR factorizations,
\begin{align*}
{\bm Q}{\bm R} & =
{\bm A} =
{\bm U}{\bm \Sigma}{\bm V}^\top \\
{\bm R} &= {\bm Q}^\top {\bm U}{\bm \Sigma}{\bm V}^\top.
\end{align*}
\begin{algorithm}[t]
\caption{Parallel Left Orthogonalization}\label{alg:ptt-l-orthog}
\KwData{A tensor $\bm f$ in distributed memory
TT format with cores
$([{\bm C}_1^p],[{\bm C}_2^p],\dots,[{\bm C}_d^p])$.}
\KwResult{A tensor ${\bm g}={\bm f}$ in
distributed memory TT format with
left orthogonal cores.}
\vspace{5mm}
Reserve memory for each block
${\bm D}_c^p$ with sizes of ${\bm C}_c^p$.\\
${\bm D}_1^p = {\bm C}_1^p$\\
\For{$c=1,2,\dots ,d-1$}{
	$[{\bm Q}_c^p,{\bm R}_c]
		= TSQR({\mathtt V}({\bm D}_c^p))$\\
	${\bm H}_{c+1}^p
		= {\bm R}_c{\mathtt H}({\bm C}_{c+1}^p)$\\
	${\bm D}_{c}^p[:] = {\bm Q}_c^p[:]$\\  
	${\bm D}_{c+1}^p[:] = {\bm H}_{c+1}^p[:]$\\  
}
Set the cores of $\bm g$ as 
$([{\bm D}_1^p],[{\bm D}_2^p],\dots,[{\bm D}_d^p])$
\end{algorithm}
Thus ${\bm R}$ has the same singular values as ${\bm A}$, and
so truncating an SVD of ${\bm R}$ produces the same approximation
error as truncating ${\bm A}$ directly. Since the TSQR algorithm
is designed to produce a copy of ${\bm R}$ on each compute node,
we apply a shared memory SVD redundantly on this copied ${\bm R}$
to determine the truncated tensor. This process is formalized
in Algorithms \ref{alg:ptt-l-orthog} and \ref{alg:ptt-r-trunc}.
These algorithms are designed so that the only communications
required are those that are involved in the computation
of the QR or LQ factorizations.
Each algorithm is mathematically equivalent to its serial
counterpart. The core difference is the use of a superscript
$p=1,2,\dots, P$,
to denote which matrices take on different values within
the different compute nodes. Matrices which are numerically
identical across all nodes lack this superscript.
\begin{algorithm}[t]
\caption{Parallel TT Truncation}\label{alg:ptt-r-trunc}
\KwData{A tensor $\bm f$ in distributed memory
TT format with cores
$([{\bm C}_1^p],[{\bm C}_2^p],\dots,[{\bm C}_d^p])$
and a desired accuracy $\varepsilon$.}
\KwResult{A tensor ${\bm g}$ so that
$\|{\bm g}-{\bm f}\|\leq \varepsilon$
in TT format with
right orthogonal cores.}
\vspace{5mm}
Reserve memory for each block
${\bm D}_c^p$ with sizes of ${\bm C}_c^p$.\\
Set $([{\hat {\bm C}}_1^p],[{\hat {\bm C}}_2^p],
\dots,[{\hat {\bm C}}_d^p])$ as orthogonal by Alg.
\ref{alg:ptt-l-orthog}.\\
${\hat \varepsilon} = \varepsilon/\sqrt{d-1}$\\
${\bm D}_d^p = {\hat {\bm C}}_d^p$\\
\For{$c=d,d-1,\dots ,2$}{
	$[{\bm L}_c,{\bm Q}_c^p] = WFLQ({\mathtt H}({\bm D}_c^p))$\\
	$[{\bm U}_c,{\bm \Sigma}_c,{\bm V}_c^\top] =
			SVD({\bm L}_c,{\hat \varepsilon})$
			\quad (Truncated SVD)\\
	${\bm H}_{c-1}^p = {\mathtt V}({\hat {\bm C}}_{c-1}^p)
				{\bm U}_c {\bm \Sigma}_c $\\
	${\bm {\hat Q}}_c^p = 	{\bm V}_c^\top{\bm Q}_c^p$\\
	${\bm D}_{c}^p[:] = {\bm {\hat Q}}_c^p [:]$\\  
	${\bm D}_{c-1}^p[:] = {\bm H}_{c-1}^p[:]$\\  
}
Set the cores of $\bm g$ as 
$([{\bm D}_1^p],[{\bm D}_2^p],\dots,[{\bm D}_d^p])$
\end{algorithm}

\bibliography{fde_tensors}

\end{document}